\documentclass[a4paper, 12pt]{article}%
\usepackage{amsmath}
\usepackage{graphicx}%
\usepackage{amsfonts}%
\usepackage{amssymb}

\begin{document}

\begin{center}
{\Huge A van der Corput lemma and weak mixing over groups}

{\Huge \vspace{10mm} }{\normalsize Conrad Beyers, Rocco
Duvenhage\footnote{{\normalsize Corresponding author.  }
\par
{\normalsize \textit{E-mail adresses: }conrad.beyers@up.ac.za (C. Beyers),
rocco@postino.up.ac.za (R. Duvenhage), anton.stroh@up.ac.za (A. Str\"{o}h).  }
\par
{\normalsize {}}} and Anton Str\"{o}h }

\emph{Department of Mathematics and Applied Mathematics}

\emph{University of Pretoria, 0002 Pretoria, South Africa}{\normalsize \ }

\bigskip

2005-8-17

\end{center}

\bigskip

\noindent\textbf{Abstract.} We study weak mixing of all orders for weakly
mixing measure preserving dynamical systems, where the dynamics is given by
the action of an abelian second countable topological group which has an
invariant measure under the group operation. One of the main technical tools
we use is a van der Corput lemma for Hilbert space valued functions on a
second countable topological group.

\newpage

\section{Introduction}

Weak mixing is an important notion in ergodic theory, introduced by Koopman
and von Neumann [9] in 1932 for actions of the group $\mathbb{R}$. Furstenberg
[5, 6] studied weakly mixing $\mathbb{Z}$-actions in order to give an ergodic
theoretic proof of Szemer\'{e}di's theorem in combinatorial number theory, and
in the process he proved that weakly mixing systems are weakly mixing of all
orders. In the case of $\mathbb{Z}$, a measure preserving transformation $T$
of a probability space $(X,\Sigma,\nu)$, namely a set $X$ with $\sigma
$-algebra $\Sigma$ on which a measure $\nu$ with $\nu(X)=1$ is defined, is
called weakly mixing if
\begin{equation}
\lim_{N\rightarrow\infty}\frac{1}{N}\sum_{n=1}^{N}\left|  \nu(A\cap
T^{-n}(B))-\nu(A)\nu(B)\right|  =0 \tag{1.1}%
\end{equation}
for all $A,B\in\Sigma$. We call this system weakly mixing of all orders if
\begin{equation}
\lim_{N\rightarrow\infty}\frac{1}{N}\sum_{n=1}^{N}\left|  \nu(A_{0}\cap
T^{-m_{1}n}(A_{1})\cap...\cap T^{-m_{k}n}(A_{k}))-\nu(A_{0})\nu(A_{1}%
)...\nu(A_{k})\right|  =0 \tag{1.2}%
\end{equation}
for all $A_{0},...,A_{k}\in\Sigma$, all $m_{1},...,m_{k}\in\mathbb{N}$ with
$m_{1}<m_{2}<...<m_{k}$, and all $k\in\mathbb{N=}\left\{  1,2,3,...\right\}  $.

However, weak mixing has also been studied for more general group (and
semigroup) actions, notably in [4] and [2] (but also see references therein),
in terms of invariant means [11] on certain spaces of functions, instead of in
terms of the explicit form $\frac{1}{N}\sum_{n=1}^{N}$ above. In particular,
various characterizations of weak mixing over the groups $\mathbb{Z}$ and
$\mathbb{R}$ were extended to more general groups.

In this paper we study weak mixing of all orders for more general group
actions than $\mathbb{Z}$ and $\mathbb{R}$. One of the technical tools we use,
is a so-called van der Corput lemma which we discuss in Section 2. This type
of lemma and related inequalities, inspired by the classical van der Corput
difference theorem and van der Corput inequality, have been used by Bergelson
[1], Furstenberg [7], Niculescu, Str\"{o}h, and Zsid\'{o} [10], and others, to
study polynomial ergodic theorems, nonconventional ergodic averages, and
noncommutative recurrence, for example. In Section 2 we extend the van der
Corput lemma to groups more general than $\mathbb{Z}$. The main results of
this section are given by Theorems 2.7 and 2.7$^{\prime}$. Instead of working
with an invariant mean on spaces of functions, we generalize the $\frac{1}%
{N}\sum_{n=1}^{N}$ form more directly for groups with an invariant measure,
since this seems convenient in our proof of the van der Corput lemma. Because
of this, we also study weak mixing by generalizing the $\frac{1}{N}\sum
_{n=1}^{N}$ form directly, rather than using the invariant mean approach of
[2]. The groups over which we work, need to have an invariant measure, and a
space-filling sequence (defined in Section 2). After some preliminaries on
weak mixing in Section 3, we devote Section 4 to showing how weak mixing
implies weak mixing of all orders, for actions of abelian second countable
topological groups of this type. The form of weak mixing of all orders we
prove, involves replacing the multiplication with $m_{1},...,m_{k}$ in (1.2),
by homomorphisms of the group over which we work. The main result is Theorem 4.4.

\section{A van der Corput lemma}

This section is devoted to proving a van der Corput lemma, stated in two
versions in Theorems 2.7 and 2.7$^{\prime}$. Our proof of the van der Corput
lemma will roughly follow that of [7] over the group $\mathbb{Z}$. In this
section we will work over a second countable topological group (i.e. a second
countable topological space which is also a group with continuous product and
inverse), since for second countable topological spaces $X,Y$, and their Borel
$\sigma$-algebras $S,T$, the product $\sigma$-algebra obtained from $S,T$ is
the same as the Borel $\sigma$-algebra of the topological space $X\times Y$.
This is needed in order to apply Fubini's theorem, which requires
measurability in the product $\sigma$-algebra. The groups that we will
consider in this paper, will only be required to be abelian from Definition
4.2 and onwards, in Section 4.

For $(Y,\mu)$ a measure space and $\mathfrak{H}$ a Hilbert space, consider a
bounded $f:\Lambda\rightarrow\mathfrak{H}$ with $\Lambda\subset Y$ measurable
and $\mu(\Lambda)<\infty$, and $\left\langle f(\cdot),x\right\rangle $
measurable for every $x\in\mathfrak{H}$. Define $\int_{\Lambda}fd\mu$ by
requiring
\[
\left\langle \int_{\Lambda}fd\mu,x\right\rangle :=\int_{\Lambda}\left\langle
f(y),x\right\rangle d\mu(y)
\]
for all $x\in\mathfrak{H}$. We will often use the notation $\int_{\Lambda
}f(y)dy=\int_{\Lambda}fd\mu$, since there will be no ambiguity in the measure
being used. Iterated integrals (when they exist) will be written as $\int
_{B}\int_{A}f(y,z)dydz$, which of course simply means $\int_{B}\left[
\int_{A}f(y,z)dy\right]  dz$, and similarly for triple integrals.

For a group $G$ we call $K\subset G$ a subsemigroup if $ab\in K$ for all
$a,b\in K$. A \textit{right invariant measure} on a topological group $G$, is
a positive measure $\mu$ on the Borel $\sigma$-algebra of $G$, with
$\mu(\Lambda g)=\mu(\Lambda)$ for all Borel $\Lambda\subset G$ and $g\in G$.
Similarly for a left invariant measure. If the measure is both right and left
invariant, we simply call it \textit{invariant}. We define such measures for
topological semigroups in the same way.

When we say that a net $\left\{  \Lambda_{\alpha}\right\}  $ has some property
for $\alpha$ ``large enough'', then we mean that there is a $\beta$ in the
directed set such that the property holds for all $\alpha\geq\beta$.

\bigskip

\noindent\textbf{Definition 2.1.} Consider a Borel measurable subsemigroup $K
$ of a topological group $G$ with right invariant measure $\mu$. A net
$\left\{  \Lambda_{\alpha}\right\}  $ of Borel subsets of $K$ is called a
\textit{space-filling net in} $K$ (or a \textit{F\o lner net in} $K$) if
$\mu(\Lambda_{\alpha})<\infty$, $\mu(\Lambda_{\beta})>0$ for $\beta$ large
enough, and
\[
\lim_{\beta}\frac{1}{\mu(\Lambda_{\beta})}\mu\left(  \Lambda_{\beta}%
\Delta(\Lambda_{\beta}g)\right)  =0
\]
for all $g\in K$. This net $\left\{  \Lambda_{\alpha}\right\}  $ is called
\textit{uniformly }space-filling if in addition%
\[
\lim_{\beta}\frac{1}{\mu(\Lambda_{\beta})}\sup_{g\in\Lambda_{\alpha}}%
\mu\left(  \Lambda_{\beta}\Delta(\Lambda_{\beta}g)\right)  =0
\]
for all $\alpha$ in the directed set of the net.

\bigskip

At the end of Section 4, we briefly consider simple examples of such nets.

\bigskip

\noindent\textbf{Proposition 2.2.} \emph{Let $G$ be a second countable
topological group with right invariant measure $\mu$. Let $\left\{
\Lambda_{\alpha}\right\}  $ be a uniformly space-filling net in a Borel
measurable subsemigroup $K$ of $G$. Consider a bounded $f:K\rightarrow
\mathfrak{H}$ with $\mathfrak{H}$ a Hilbert space, such that $\left\langle
f(\cdot),x\right\rangle $ is Borel measurable for every $x\in\mathfrak{H}$.
Then}
\[
\lim_{\beta}\left|  \left|  \frac{1}{\mu(\Lambda_{\beta})}\int_{\Lambda
_{\beta}}fd\mu-\frac{1}{\mu(\Lambda_{\beta})}\frac{1}{\mu(\Lambda_{\alpha}%
)}\int_{\Lambda_{\beta}}\int_{\Lambda_{\alpha}}f(gh)dhdg\right|  \right|  =0
\]
\emph{for every $\alpha$ in the directed set of the net.}

\bigskip

\noindent\textbf{Proof.} By Fubini's theorem
\begin{align*}
\int_{\Lambda_{\beta}}\int_{\Lambda_{\alpha}}\left\langle f(gh),x\right\rangle
dhdg  &  =\int_{\Lambda_{\beta}\times\Lambda_{\alpha}}\left\langle
f(gh),x\right\rangle d(g,h)\\
&  =\int_{\Lambda_{\alpha}}\int_{\Lambda_{\beta}}\left\langle
f(gh),x\right\rangle dgdh
\end{align*}
which by definition means that%
\[
\int_{\Lambda_{\beta}}\int_{\Lambda_{\alpha}}f(gh)dhdg=\int_{\Lambda_{\alpha}%
}\int_{\Lambda_{\beta}}f(gh)dgdh
\]
and in particular these iterated integrals exists. From this and the fact that
$\mu$ is a right invariant measure, we have
\begin{align*}
&  \left|  \left|  \frac{1}{\mu(\Lambda_{\beta})}\int_{\Lambda_{\beta}}%
fd\mu-\frac{1}{\mu(\Lambda_{\beta})}\frac{1}{\mu(\Lambda_{\alpha})}%
\int_{\Lambda_{\beta}}\int_{\Lambda_{\alpha}}f(gh)dhdg\right|  \right| \\
&  =\left|  \left|  \frac{1}{\mu(\Lambda_{\alpha})}\int_{\Lambda_{\alpha}%
}\left[  \frac{1}{\mu(\Lambda_{\beta})}\int_{\Lambda_{\beta}}f(g)dg\right]
dh-\frac{1}{\mu(\Lambda_{\alpha})}\frac{1}{\mu(\Lambda_{\beta})}\int
_{\Lambda_{\alpha}}\left[  \int_{\Lambda_{\beta}}f(gh)dg\right]  dh\right|
\right| \\
&  =\left|  \left|  \frac{1}{\mu(\Lambda_{\alpha})}\frac{1}{\mu(\Lambda
_{\beta})}\int_{\Lambda_{\alpha}}\left[  \int_{\Lambda_{\beta}}f(g)dg-\int
_{\Lambda_{\beta}}f(gh)dg\right]  dh\right|  \right| \\
&  =\left|  \left|  \frac{1}{\mu(\Lambda_{\alpha})}\frac{1}{\mu(\Lambda
_{\beta})}\int_{\Lambda_{\alpha}}\left[  \int_{\Lambda_{\beta}}f(g)dg-\int
_{\Lambda_{\beta}h}f(g)dg\right]  dh\right|  \right| \\
&  =\left|  \left|  \frac{1}{\mu(\Lambda_{\alpha})}\frac{1}{\mu(\Lambda
_{\beta})}\int_{\Lambda_{\alpha}}\left[  \int_{\Lambda_{\beta}\backslash
\left(  \Lambda_{\beta}\cap(\Lambda_{\beta}h)\right)  }f(g)dg-\int_{\left(
\Lambda_{\beta}h\right)  \backslash\left(  \Lambda_{\beta}\cap(\Lambda_{\beta
}h)\right)  }f(g)dg\right]  dh\right|  \right|  \text{ \ \ .}%
\end{align*}
But if$\ b\in\mathbb{R}$ is an upper bound for $\left|  \left|  f(K)\right|
\right|  $ , we have%
\begin{align*}
&  \left|  \left|  \int_{\Lambda_{\beta}\backslash\left(  \Lambda_{\beta}%
\cap(\Lambda_{\beta}h)\right)  }f(g)dg-\int_{\left(  \Lambda_{\beta}h\right)
\backslash\left(  \Lambda_{\beta}\cap(\Lambda_{\beta}h)\right)  }%
f(g)dg\right|  \right| \\
&  \leq b\mu\left(  \Lambda_{\beta}\backslash\left(  \Lambda_{\beta}%
\cap(\Lambda_{\beta}h)\right)  \right)  +b\mu\left(  \left(  \Lambda_{\beta
}h\right)  \backslash\left(  \Lambda_{\beta}\cap(\Lambda_{\beta}h)\right)
\right) \\
&  =b\mu\left(  \Lambda_{\beta}\Delta(\Lambda_{\beta}h)\right) \\
&  \leq b\sup_{h\in\Lambda_{\alpha}}\mu\left(  \Lambda_{\beta}\Delta
(\Lambda_{\beta}h)\right)
\end{align*}
therefore%
\begin{align*}
&  \left|  \left|  \frac{1}{\mu(\Lambda_{\beta})}\int_{\Lambda_{\beta}}%
fd\mu-\frac{1}{\mu(\Lambda_{\beta})}\frac{1}{\mu(\Lambda_{\alpha})}%
\int_{\Lambda_{\beta}}\int_{\Lambda_{\alpha}}f(gh)dhdg\right|  \right| \\
&  \leq\frac{1}{\mu(\Lambda_{\beta})}b\sup_{h\in\Lambda_{\alpha}}\mu\left(
\Lambda_{\beta}\Delta(\Lambda_{\beta}h)\right) \\
&  \rightarrow0
\end{align*}
in the $\beta$ limit.$~\square$

\bigskip

\noindent\textbf{Lemma 2.3.}\emph{\ Let $\mathfrak{H}$ be a Hilbert space,
$(Y,\mu)$ a measure space, and }$\Lambda\subset Y$\emph{\ a measurable set
with $\mu(\Lambda)<\infty$. Consider an $f:\Lambda\rightarrow\mathfrak{H}$
with $\left|  \left|  f(\cdot)\right|  \right|  $ measurable, and
$\left\langle f(\cdot),x\right\rangle $ measurable for every $x\in
\mathfrak{H}$, and with $\int_{\Lambda}\left|  \left|  f(y)\right|  \right|
dy<\infty$ (which means $\int_{\Lambda}fd\mu$ exists). Then}%
\[
\left|  \left|  \int_{\Lambda}fd\mu\right|  \right|  ^{2}\leq\mu(\Lambda
)\int_{\Lambda}\left|  \left|  f(y)\right|  \right|  ^{2}dy\text{ \ \ .}%
\]

\noindent\textbf{Proof.} By definition of $\int_{\Lambda}fd\mu$,
\begin{align*}
\left|  \left|  \int_{\Lambda}fd\mu\right|  \right|  ^{2}  &  =\left\langle
\int_{\Lambda}fd\mu,\int_{\Lambda}fd\mu\right\rangle =\int_{\Lambda
}\left\langle f(y),\int_{\Lambda}fd\mu\right\rangle dy\\
&  =\int_{\Lambda}\left[  \int_{\Lambda}\left\langle f(y),f(z)\right\rangle
dz\right]  dy\text{ \ \ .}%
\end{align*}
For any $a,b\in\mathfrak{H}$ we have $2\operatorname{Re}\left\langle
a,b\right\rangle \leq\left|  \left|  a\right|  \right|  ^{2}+\left|  \left|
b\right|  \right|  ^{2}$, and since the object above is real, we have%
\begin{align*}
\left|  \left|  \int_{\Lambda}fd\mu\right|  \right|  ^{2}  &  =\int_{\Lambda}
\left[  \int_{\Lambda}\operatorname{Re}\left\langle f(y),f(z)\right\rangle
dz\right]  dy\\
&  \leq\frac{1}{2}\int_{\Lambda}\left[  \int_{\Lambda}\left(  \left|  \left|
f(y)\right|  \right|  ^{2}+\left|  \left|  f(z)\right|  \right|  ^{2}\right)
dz\right]  dy\\
&  =\mu(\Lambda)\int_{\Lambda}\left|  \left|  f(y)\right|  \right|
^{2}dy\text{ \ \ .}~\square
\end{align*}

\bigskip

\noindent\textbf{Proposition 2.4.}\emph{\ Consider the situation in
Proposition 2.2, except that we don't need the net. Assume furthermore that
$F:K\times K\rightarrow\mathbb{C}:(g,h)\mapsto\left\langle
f(g),f(h)\right\rangle $ is Borel measurable, and that }$\Lambda_{1}%
,\Lambda_{2}\subset K$\emph{\ are Borel sets with }$\mu(\Lambda_{j})<\infty
$\emph{. Then}%
\begin{align*}
&  \left|  \left|  \int_{\Lambda_{2}}\int_{\Lambda_{1}}f(gh)dhdg\right|
\right|  ^{2}\\
&  \leq\mu(\Lambda_{2})\int_{\Lambda_{1}}\int_{\Lambda_{1}}\int_{\Lambda_{2}%
}\left\langle f(gh_{1}),f(gh_{2})\right\rangle dgdh_{1}dh_{2}%
\end{align*}
\emph{and in particular these integrals exist.}

\bigskip

\noindent\textbf{Proof.} The double integral exists as in Proposition 2.2's
proof. Let's now consider the triple integral. Since $F$ is Borel measurable
and $G$ 's product is continuous, $(g,h_{1})\mapsto\left\langle f(gh_{1}%
),f(gh_{2})\right\rangle $ is Borel measurable on $K\times K=K^{2}$ and hence
measurable in the product $\sigma$-algebra on $K^{2}$. By Fubini's theorem we
have%
\[
\int_{\Lambda_{1}}\int_{\Lambda_{2}}\left\langle f(gh_{1}),f(gh_{2}%
)\right\rangle dgdh_{1}=\int_{\Lambda_{1}\times\Lambda_{2}}\left\langle
f(gh_{1}),f(gh_{2})\right\rangle d(h_{1},g)
\]
and in particular the iterated integral exists. Furthermore, $K\times
K^{2}\rightarrow K^{2}:(h_{2},h_{1},g)\mapsto(gh_{1},gh_{2})$ is continuous,
so $K\times K^{2}\rightarrow\mathbb{C}:(h_{2},h_{1},g)\mapsto\left\langle
f(gh_{1}),f(gh_{2})\right\rangle $ is measurable in the product $\sigma
$-algebra of $K$ and $K^{2}$. Hence by Fubini's theorem
\[
\int_{\Lambda_{1}}\int_{\Lambda_{1}\times\Lambda_{2}}\left\langle
f(gh_{1}),f(gh_{2})\right\rangle d(h_{1},g)dh_{2}=\int_{\Lambda_{1}%
\times\Lambda_{1}\times\Lambda_{2}}\left\langle f(gh_{1}),f(gh_{2}%
)\right\rangle d(h_{2},h_{1},g)
\]
and in particular, the triple integral exists, and we can do the three
integrals in any order. By Lemma 2.3 it follows that
\begin{align*}
&  \left|  \left|  \int_{\Lambda_{2}}\int_{\Lambda_{1}}f(gh)dhdg\right|
\right|  ^{2}\\
&  \leq\mu(\Lambda_{2})\int_{\Lambda_{2}}\left|  \left|  \int_{\Lambda_{1}%
}f(gh)dh\right|  \right|  ^{2}dg\\
&  =\mu(\Lambda_{2})\int_{\Lambda_{2}}\left\langle \int_{\Lambda_{1}}%
f(gh_{1})dh_{1},\int_{\Lambda_{1}}f(gh_{2})dh_{2}\right\rangle dg\\
&  =\mu(\Lambda_{2})\int_{\Lambda_{2}}\int_{\Lambda_{1}}\left\langle
f(gh_{1}),\int_{\Lambda_{1}}f(gh_{2})dh_{2}\right\rangle dh_{1}dg\\
&  =\mu(\Lambda_{2})\int_{\Lambda_{2}}\int_{\Lambda_{1}}\int_{\Lambda_{1}%
}\left\langle f(gh_{1}),f(gh_{2})\right\rangle dh_{2}dh_{1}dg\\
&  =\mu(\Lambda_{2})\int_{\Lambda_{1}}\int_{\Lambda_{1}}\int_{\Lambda_{2}%
}\left\langle f(gh_{1}),f(gh_{2})\right\rangle dgdh_{1}dh_{2}%
\end{align*}
and note in particular that the part of this argument after the inequality
proves that $g\mapsto\left|  \left|  \int_{\Lambda_{1}}f(gh)dh\right|
\right|  ^{2}$ is measurable (and therefore its square root too), which means
that Lemma 2.3 does indeed apply to this situation.$~\square$

\bigskip

For the next three results we give two versions of each, one set of results
for nets, and one for sequences but with other assumptions a bit weaker. The
weaker assumptions in case of sequences are possible, since in this case we
can apply Lebesgue's dominated convergence theorem (see the proof of
Proposition 2.6$^{\prime}$). The case of sequences will be used in Section 4
when we study weak mixing of all orders.

\bigskip

\noindent\textbf{Lemma 2.5.}\emph{\ Consider the situation in Proposition 2.2,
but assume that $f:G\rightarrow\mathfrak{H}$ is bounded and $\left\langle
f(\cdot),x\right\rangle $ measurable for all $x\in\mathfrak{H}$. The net is
still uniformly space-filling only in }$\emph{K}$\emph{, though. Assume that
$F:G^{2}\rightarrow\mathbb{C}:(g,h)\mapsto\left\langle f(g),f(h)\right\rangle
$ is Borel measurable. Then $\int_{\Lambda}\left\langle
f(g),f(gh)\right\rangle dg$ exists for all measurable $\Lambda\subset G$ with
$\mu(\Lambda)<\infty$, and all $h\in G$. Assume that}%
\[
\gamma_{h}:=\lim_{\beta}\frac{1}{\mu(\Lambda_{\beta})}\int_{\Lambda_{\beta}%
}\left\langle f(g),f(gh)\right\rangle dg
\]
\emph{exists for all $h\in G$, and that}
\begin{equation}
\lim_{\beta}\sup_{h\in\Lambda_{\alpha}}\left|  \frac{1}{\mu(\Lambda_{\beta}%
)}\int_{\Lambda_{\beta}}\left\langle f(g),f(gh)\right\rangle dg-\gamma
_{h}\right|  =0 \tag{2.5.1}%
\end{equation}
\emph{for all $\alpha$. Then}%
\begin{equation}
\lim_{\beta}\sup_{h_{1},h_{2}\in\Lambda_{\alpha}}\left|  \frac{1}{\mu
(\Lambda_{\beta})}\int_{\Lambda_{\beta}}\left\langle f(gh_{1}),f(gh_{2}%
)\right\rangle dg-\gamma_{h_{1}^{-1}h_{2}}\right|  =0 \tag{2.5.2}%
\end{equation}
\emph{for all $\alpha$. In particular $\Lambda_{\alpha}\times\Lambda_{\alpha
}\ni(h_{1},h_{2})\mapsto\gamma_{h_{1}^{-1}h_{2}}$ is bounded for every
$\alpha$.}

\bigskip

\noindent\textbf{Proof.} Since $F$ is Borel, and $G\rightarrow G^{2}%
:g\mapsto(g,gh)$ is continuous, the map $G\rightarrow\mathbb{C}:g\mapsto
\left\langle f(g),f(gh)\right\rangle $ is Borel for every $h\in G$, hence the
integrals are defined. Now,%
\begin{align*}
&  \sup_{h_{1},h_{2}\in\Lambda_{\alpha}}\left|  \frac{1}{\mu(\Lambda_{\beta}%
)}\int_{\Lambda_{\beta}}\left\langle f(gh_{1}),f(gh_{2})\right\rangle
dg-\gamma_{h_{1}^{-1}h_{2}}\right| \\
&  \leq\sup_{h_{1},h_{2}\in\Lambda_{\alpha}}\frac{1}{\mu(\Lambda_{\beta}%
)}\left|  \int_{\Lambda_{\beta}}\left\langle f(gh_{1}),f(gh_{2})\right\rangle
dg-\int_{\Lambda_{\beta}}\left\langle f(g),f(gh_{1}^{-1}h_{2})\right\rangle
dg\right| \\
&  +\sup_{h_{1},h_{2}\in\Lambda_{\alpha}}\left|  \frac{1}{\mu(\Lambda_{\beta
})}\int_{\Lambda_{\beta}}\left\langle f(g),f(gh_{1}^{-1}h_{2})\right\rangle
dg-\gamma_{h_{1}^{-1}h_{2}}\right|
\end{align*}
but since $\mu$ is a right invariant measure%
\begin{align*}
&  \frac{1}{\mu(\Lambda_{\beta})}\left|  \int_{\Lambda_{\beta}}\left\langle
f(gh_{1}),f(gh_{2})\right\rangle dg-\int_{\Lambda_{\beta}}\left\langle
f(g),f(gh_{1}^{-1}h_{2})\right\rangle dg\right| \\
&  =\frac{1}{\mu(\Lambda_{\beta})}\left|  \int_{\Lambda_{\beta}h_{1}%
}\left\langle f(g),f(gh_{1}^{-1}h_{2})\right\rangle dg-\int_{\Lambda_{\beta}%
}\left\langle f(g),f(gh_{1}^{-1}h_{2})\right\rangle dg\right| \\
&  =\frac{1}{\mu(\Lambda_{\beta})}\left|  \int_{(\Lambda_{\beta}%
h_{1})\backslash\Lambda_{\beta}}\left\langle f(g),f(gh_{1}^{-1}h_{2}%
)\right\rangle dg-\int_{\Lambda_{\beta}\backslash(\Lambda_{\beta}h_{1}%
)}\left\langle f(g),f(gh_{1}^{-1}h_{2})\right\rangle dg\right| \\
&  \leq\frac{1}{\mu(\Lambda_{\beta})}\left[  \int_{(\Lambda_{\beta}%
h_{1})\backslash\Lambda_{\beta}}\left|  \left\langle f(g),f(gh_{1}^{-1}%
h_{2})\right\rangle \right|  dg+\int_{\Lambda_{\beta}\backslash(\Lambda
_{\beta}h_{1})}\left|  \left\langle f(g),f(gh_{1}^{-1}h_{2})\right\rangle
\right|  dg\right] \\
&  =\frac{1}{\mu(\Lambda_{\beta})}\int_{\Lambda_{\beta}\Delta(\Lambda_{\beta
}h_{1})}\left|  \left\langle f(g),f(gh_{1}^{-1}h_{2})\right\rangle \right|
dg\\
&  \leq\frac{\mu\left(  \Lambda_{\beta}\Delta(\Lambda_{\beta}h_{1})\right)
}{\mu(\Lambda_{\beta})}b
\end{align*}
for all $h_{1}\in K$, where $b$ is an upper bound for $(g,h_{1},h_{2}%
)\mapsto\left|  \left\langle f(g),f(gh_{1}^{-1}h_{2})\right\rangle \right|  $,
which exists since $f$ is bounded. This proves (2.5.2). Since $f$ is bounded,
so is $(h_{1},h_{2})\mapsto\left\langle f(gh_{1}),f(gh_{2})\right\rangle $ and
its integral with respect to $g$ over $\Lambda_{\beta}$. Hence (2.5.2) implies
that $\Lambda_{\alpha}\times\Lambda_{\alpha}\ni(h_{1},h_{2})\mapsto
\gamma_{h_{1}^{-1}h_{2}}$ is bounded.$~\square$

\bigskip

\noindent\textbf{Lemma 2.5}$^{\prime}$\textbf{.} \emph{Consider the situation
in Lemma 2.5, except that $\left\{  \Lambda_{\alpha}\right\}  $ need not be
uniform. Assume}%
\[
\gamma_{h}:=\lim_{\beta}\frac{1}{\mu(\Lambda_{\beta})}\int_{\Lambda_{\beta}%
}\left\langle f(g),f(gh)\right\rangle dg
\]
\emph{exists for all $h\in G$. (We need not assume 2.5.1.) Then}%
\[
\lim_{\beta}\frac{1}{\mu(\Lambda_{\beta})}\int_{\Lambda_{\beta}}\left\langle
f(gh_{1}),f(gh_{2})\right\rangle dg=\gamma_{h_{1}^{-1}h_{2}}%
\]
\emph{for all $h_{1}\in K$ and $h_{2}\in G$.}

\bigskip

\noindent\textbf{Proof.} Simply repeat Lemma 2.5 's proof without the $\sup$
's.$~\square$

\bigskip

\noindent\textbf{Proposition 2.6.} \emph{Consider the situation in Lemma 2.5.
Assuming that $K^{2}\rightarrow\mathbb{C}:(h_{1},h_{2})\mapsto\gamma
_{h_{1}^{-1}h_{2}}$ is Borel measurable, we have}
\[
\lim_{\beta}\frac{1}{\mu(\Lambda_{\beta})}\int_{\Lambda_{\alpha}}\int
_{\Lambda_{\alpha}}\int_{\Lambda_{\beta}}\left\langle f(gh_{1}),f(gh_{2}%
)\right\rangle dgdh_{1}dh_{2}=\int_{\Lambda_{\alpha}}\int_{\Lambda_{\alpha}%
}\gamma_{h_{1}^{-1}h_{2}}dh_{1}dh_{2}%
\]
\emph{for all $\alpha$.}

\bigskip

\noindent\textbf{Proof.} The triple integral exists by Proposition 2.4. The
double integral exists by Fubini's theorem, since $(h_{1},h_{2})\mapsto
\gamma_{h_{1}^{-1}h_{2}}$ is bounded on $\Lambda_{\alpha}\times\Lambda
_{\alpha}$ for every $\alpha$ by Lemma 2.5. Then%
\begin{align*}
&  \left|  \frac{1}{\mu(\Lambda_{\beta})}\int_{\Lambda_{\alpha}}\int
_{\Lambda_{\alpha}}\int_{\Lambda_{\beta}}\left\langle f(gh_{1}),f(gh_{2}%
)\right\rangle dgdh_{1}dh_{2}-\int_{\Lambda_{\alpha}}\int_{\Lambda_{\alpha}%
}\gamma_{h_{1}^{-1}h_{2}}dh_{1}dh_{2}\right| \\
&  \leq\mu(\Lambda_{\alpha})^{2}\sup_{h_{1},h_{2}\in\Lambda_{\alpha}}\left|
\frac{1}{\mu(\Lambda_{\beta})}\int_{\Lambda_{\beta}}\left\langle
f(gh_{1}),f(gh_{2})\right\rangle dgdh_{1}dh_{2}-\gamma_{h_{1}^{-1}h_{2}%
}\right| \\
&  \rightarrow0
\end{align*}
in the $\beta$ limit by Lemma 2.5.$~\square$

\bigskip

\noindent\textbf{Proposition 2.6}$^{\prime}$\textbf{.} \emph{Consider the
situation in Lemma 2.5$^{\prime}$, but assume the space-filling net in }%
$K$\emph{\ is in fact a sequence $\left\{  \Lambda_{n}\right\}  _{n\in
\mathbb{N}} $. Then}
\[
\lim_{n\rightarrow\infty}\frac{1}{\mu(\Lambda_{n})}\int_{\Lambda_{m}}%
\int_{\Lambda_{m}}\int_{\Lambda_{n}}\left\langle f(gh_{1}),f(gh_{2}%
)\right\rangle dgdh_{1}dh_{2}=\int_{\Lambda_{m}}\int_{\Lambda_{m}}%
\gamma_{h_{1}^{-1}h_{2}}dh_{1}dh_{2}%
\]
\emph{for all $m$, and in particular these integrals exist.}

\bigskip

\noindent\textbf{Proof.} The triple integral exists by Proposition 2.4. Let
$b$ be an upper bound for $(g,h_{1},h_{2})\mapsto\left|  \left\langle
f(gh_{1}),f(gh_{2})\right\rangle \right|  $, which exists since $f$ is
bounded. Fix any $m\in\mathbb{N}$, and set
\[
A_{n}(h_{1},h_{2}):=\frac{1}{\mu(\Lambda_{n})}\int_{\Lambda_{n}}\left\langle
f(gh_{1}),f(gh_{2})\right\rangle dg
\]
for all $h_{1},h_{2}\in\Lambda_{m}$ and all $n$. Note that $A_{n}(h_{1}%
,h_{2})$ exists and is a measurable function of $h_{1}$ because of the
existence of the triple integral. Then $\left|  A_{n}(h_{1},h_{2})\right|
\leq\frac{1}{\mu(\Lambda_{n})}\int_{\Lambda_{n}}bdg=b$, which implies that the
sequence $A_{n}(\cdot,h_{2})$ with $h_{2}$ fixed, is dominated by
$B:\Lambda_{m}\rightarrow\mathbb{R}:h\mapsto b$. But $B\in L^{1}(\Lambda
_{m},\mu)$, namely $\int_{\Lambda_{m}}\left|  B\right|  d\mu=b\mu(\Lambda
_{m})<\infty$, hence $\Lambda_{m}\ni h_{1}\mapsto\gamma_{h_{1}^{-1}h_{2}}$ is
in $L^{1}(\Lambda_{m},\mu)$ and
\[
\lim_{n\rightarrow\infty}\int_{\Lambda_{m}}A_{n}(h_{1},h_{2})dh_{1}%
=\int_{\Lambda_{m}}\gamma_{h_{1}^{-1}h_{2}}dh_{1}%
\]
by Lebesgue's dominated convergence theorem and Lemma 2.5$^{\prime}$. Note in
particular that the last integral exists. Now set
\[
C_{n}(h_{2}):=\int_{\Lambda_{m}}A_{n}(h_{1},h_{2})dh_{1}%
\]
for all $h_{2}\in\Lambda_{m}$, and keep in mind that this exists and is a
measurable function of $h_{2}$ by Proposition 2.4, as for $A_{n}(\cdot,h_{2})
$ earlier. Then $\left|  C_{n}(h_{2})\right|  \leq\int_{\Lambda_{m}}%
bdh_{1}\leq\mu(\Lambda_{m})b$, so the sequence $C_{n}$ is dominated by
$D:\Lambda_{m}\rightarrow\mathbb{R}:h\mapsto\mu(\Lambda_{m})b$, and $D\in
L^{1}(\Lambda_{m},\mu)$. Hence by Lebesgue's dominated convergence theorem,
the function
\[
\Lambda_{m}\ni h_{2}\mapsto\int_{\Lambda_{m}}\gamma_{h_{1}^{-1}h_{2}}dh_{1}%
\]
is in $L^{1}(\Lambda_{m},\mu)$, and
\[
\lim_{n\rightarrow\infty}\int_{\Lambda_{m}}C_{n}(h_{2})dh_{2}=\int
_{\Lambda_{m}}\int_{\Lambda_{m}}\gamma_{h_{1}^{-1}h_{2}}dh_{1}dh_{2}%
\]
as required. $\square$

\bigskip

Now we can finally state a van der Corput lemma:

\bigskip

\noindent\textbf{Theorem 2.7.}\emph{\ Let $G$ be a second countable
topological group with right invariant measure $\mu$. Let $\left\{
\Lambda_{\alpha}\right\}  $ be a uniformly space-filling net in a Borel
measurable subsemigroup $K$ of $G$. Consider a bounded $f:G\rightarrow
\mathfrak{H}$, with $\mathfrak{H}$ a Hilbert space, such that $\left\langle
f(\cdot),x\right\rangle $ and $\left\langle f(\cdot),f(\cdot)\right\rangle
:G^{2}\rightarrow\mathbb{C}$ are Borel measurable (for all $x\in\mathfrak{H}%
$). Assume that}%
\[
\gamma_{h}:=\lim_{\beta}\frac{1}{\mu(\Lambda_{\beta})}\int_{\Lambda_{\beta}%
}\left\langle f(g),f(gh)\right\rangle dg
\]
\emph{exists for all $h\in G$, and that}
\[
\lim_{\beta}\sup_{h\in\Lambda_{\alpha}}\left|  \frac{1}{\mu(\Lambda_{\beta}%
)}\int_{\Lambda_{\beta}}\left\langle f(g),f(gh)\right\rangle dg-\gamma
_{h}\right|  =0
\]
\emph{for all $\alpha$. Assume furthermore $G\rightarrow\mathbb{C}%
:h\mapsto\gamma_{h}$ is Borel measurable and that}
\begin{equation}
\lim_{\alpha}\frac{1}{\mu(\Lambda_{\alpha})^{2}}\int_{\Lambda_{\alpha}}%
\int_{\Lambda_{\alpha}}\gamma_{h_{1}^{-1}h_{2}}dh_{1}dh_{2}=0\text{ \ \ .}
\tag{2.7.1}%
\end{equation}
\emph{Then}%
\[
\lim_{\beta}\frac{1}{\mu(\Lambda_{\beta})}\int_{\Lambda_{\beta}}fd\mu=0\text{
\ \ .}%
\]

\noindent\textbf{Proof.} Note that since $G^{2}\rightarrow G:(h_{1}%
,h_{2})\mapsto h_{1}^{-1}h_{2}$ is continuous, its composition with
$h\mapsto\gamma_{h}$, namely $(h_{1},h_{2})\mapsto\gamma_{h_{1}^{-1}h_{2}}$,
is Borel. By Proposition 2.2 and Proposition 2.4 we just have to show that for
any $\varepsilon>0$ there is an $\alpha$ and $\beta_{0}$ such that $\left|
A_{\alpha\beta}\right|  <\varepsilon$ for all $\beta>\beta_{0}$ where%
\[
A_{\alpha\beta}:=\frac{1}{\mu(\Lambda_{\beta})}\frac{1}{\mu(\Lambda_{\alpha
})^{2}}\int_{\Lambda_{\alpha}}\int_{\Lambda_{\alpha}}\int_{\Lambda_{\beta}%
}\left\langle f(gh_{1}),f(gh_{2})\right\rangle dgdh_{1}dh_{2}\text{ \ \ .}%
\]
But this follows from Proposition 2.6 and our assumptions, namely%
\[
\lim_{\alpha}\lim_{\beta}A_{\alpha\beta}=\lim_{\alpha}\frac{1}{\mu
(\Lambda_{\alpha})^{2}}\int_{\Lambda_{\alpha}}\int_{\Lambda_{\alpha}}%
\gamma_{h_{1}^{-1}h_{2}}dh_{1}dh_{2}=0\text{ \ \ .}~\square
\]

\bigskip

Next we give a van der Corput lemma for a space-filling sequence instead of a net:

\bigskip

\noindent\textbf{Theorem 2.7}$^{\prime}$\textbf{.} \emph{Let $G$ be a second
countable topological group with right invariant measure $\mu$. Let $\left\{
\Lambda_{n}\right\}  $ be a uniformly space-filling sequence in a Borel
measurable subsemigroup $K$ of $G$. Consider a bounded $f:G\rightarrow
\mathfrak{H}$, with $\mathfrak{H}$ a Hilbert space, such that $\left\langle
f(\cdot),x\right\rangle $ and $\left\langle f(\cdot),f(\cdot)\right\rangle
:G^{2}\rightarrow\mathbb{C}$ are Borel measurable (for all $x\in\mathfrak{H}%
$). Assume}%
\[
\gamma_{h}:=\lim_{n\rightarrow\infty}\frac{1}{\mu(\Lambda_{n})}\int
_{\Lambda_{n}}\left\langle f(g),f(gh)\right\rangle dg
\]
\emph{exists for all $h\in G$. Also assume that}
\begin{equation}
\lim_{m\rightarrow\infty}\frac{1}{\mu(\Lambda_{m})^{2}}\int_{\Lambda_{m}}%
\int_{\Lambda_{m}}\gamma_{h_{1}^{-1}h_{2}}dh_{1}dh_{2}=0 \tag{2.7$^\prime$.1}%
\end{equation}
\emph{(note that the integral exists by Proposition 2.6$^{\prime}$). Then}%
\[
\lim_{n\rightarrow\infty}\frac{1}{\mu(\Lambda_{n})}\int_{\Lambda_{n}}%
fd\mu=0\text{ \ \ .}%
\]

\noindent\textbf{Proof.} Just as for Theorem 2.7, but using Proposition
2.6$^{\prime}$ instead of 2.6, and therefore without the need to show that
$(h_{1},h_{2})\mapsto\gamma_{h_{1}^{-1}h_{2}}$ is Borel.$~\square$

\bigskip

Theorem 2.7$^{\prime}$ is the version of the van der Corput lemma that we will
apply in Section 4 to prove that weak mixing implies weak mixing of all
orders. However, we still need a few refinements regarding conditions (2.7.1)
and (2.7$^{\prime}$.1):

\bigskip

\noindent\textbf{Lemma 2.8.}\emph{\ Let $G$ be a second countable topological
group with left invariant measure $\mu$. Let $\Lambda\subset G$ be Borel and
$\mu(\Lambda)<\infty$, and $S\subset G$ Borel such that $\Lambda^{-1}%
\Lambda:=\left\{  h_{1}^{-1}h_{2}:h_{1},h_{2}\in\Lambda\right\}  \subset S$.
For a Borel $f:G\rightarrow\mathbb{R}^{+}$ we then have}%
\[
\int_{\Lambda}\int_{\Lambda}f(h_{1}^{-1}h_{2})dh_{1}dh_{2}\leq\mu(\Lambda
)\int_{S}fd\mu\text{ \ \ .}%
\]

\noindent\textbf{Proof.} Let $\chi$ denote characteristic functions, and set
$\varphi:\Lambda\times\Lambda\rightarrow G:(h_{1},h_{2})\mapsto h_{1}%
^{-1}h_{2}$. Then $f\circ\varphi$ is Borel on $\Lambda\times\Lambda$, and
therefore measurable in the product $\sigma$-algebra on $\Lambda\times\Lambda$
obtained from $\Lambda$ 's Borel $\sigma$-algebra, since $\varphi$ is
continuous. Let $Y\subset\Lambda^{-1}\Lambda$ be Borel in $G$. For $W\subset
G\times G$, let $W_{g}:=\left\{  h:(g,h)\in W\right\}  $. Then, since
$\varphi^{-1}(Y)$ is Borel in $\Lambda\times\Lambda$ and hence Borel in
$G\times G$, it follows that $\varphi^{-1}(Y)$ is in the product $\sigma
$-algebra on $G\times G$, hence we can consider $(\mu\times\mu)\left(
\varphi^{-1}(Y)\right)  =\int_{\Lambda}\mu\left(  \varphi^{-1}(Y)_{g}\right)
dg $. Now%
\[
\varphi^{-1}(Y)=\left\{  (g,gh):h\in Y,g\in\Lambda\cap\left(  \Lambda
h^{-1}\right)  \right\}  \subset\left\{  (g,gh):h\in Y,g\in\Lambda\right\}
=:V
\]
but $V_{g}=gY$, therefore $\mu\left(  \varphi^{-1}(Y)_{g}\right)  \leq
\mu(V_{g})=\mu(gY)=\mu(Y)$, since $\mu$ is a left invariant. Hence
\begin{align*}
\int_{\Lambda\times\Lambda}\chi_{Y}\circ\varphi d(\mu\times\mu)  &
=(\mu\times\mu)\left(  \varphi^{-1}(Y)\right) \\
&  \leq\mu(\Lambda)\mu(Y)\\
&  =\mu(\Lambda)\int_{S}\chi_{Y}d\mu
\end{align*}

There is an increasing sequence $f_{n}:S\rightarrow\mathbb{R}^{+}$ of simple
functions converging pointwise to $f$. From the above we know that%
\[
\int_{\Lambda\times\Lambda}f_{n}\circ\varphi d(\mu\times\mu)\leq\mu
(\Lambda)\int_{S}f_{n}d\mu
\]
and by applying Lebesgue's monotone convergence first on the right and then of
the left of this inequality, we obtain%
\[
\int_{\Lambda}\int_{\Lambda}f\left(  h_{1}^{-1}h_{2}\right)  dh_{1}dh_{2}%
=\int_{\Lambda\times\Lambda}f\circ\varphi d(\mu\times\mu)\leq\mu(\Lambda
)\int_{S}fd\mu
\]
as required, where we have used Fubini's theorem, which holds in this case,
since $f$ is non-negative.$~\square$

\bigskip

\noindent\textbf{Proposition 2.9.}\emph{\ Let $G$ be a second countable
topological group with left invariant measure $\mu$. Let $\left\{
\Lambda_{\alpha}\right\}  $ be a uniformly space-filling net in a Borel
measurable subsemigroup $K$ of $G$. Consider a Borel measurable function
$\gamma_{h}:G\rightarrow\mathbb{C}$. Also assume that each $\Lambda_{\alpha}$
is open, and that}%
\[
\lim_{\alpha}\frac{1}{\mu(\Lambda_{\alpha})}\int_{\Lambda_{\alpha}^{-1}%
\Lambda_{\alpha}}\left|  \gamma_{h}\right|  dh=0\text{ \ \ .}%
\]
\emph{Then}
\[
\lim_{\alpha}\frac{1}{\mu(\Lambda_{\alpha})^{2}}\int_{\Lambda_{\alpha}}%
\int_{\Lambda_{\alpha}}\gamma_{h_{1}^{-1}h_{2}}dh_{1}dh_{2}=0
\]
\emph{if the iterated integral exists for all $\alpha\geq\alpha_{0}$ for some
$\alpha_{0}$.}

\bigskip

\noindent\textbf{Proof.} Since $\Lambda_{\alpha}$ is open, $\Lambda_{\alpha
}^{-1}\Lambda_{\alpha}$ is Borel, and so
\begin{align*}
\left|  \frac{1}{\mu(\Lambda_{\alpha})^{2}}\int_{\Lambda_{\alpha}}%
\int_{\Lambda_{\alpha}}\gamma_{h_{1}^{-1}h_{2}}dh_{1}dh_{2}\right|   &
\leq\frac{1}{\mu(\Lambda_{\alpha})^{2}}\int_{\Lambda_{\alpha}}\int
_{\Lambda_{\alpha}}\left|  \gamma_{h_{1}^{-1}h_{2}}\right|  dh_{1}dh_{2}\\
&  \leq\frac{1}{\mu(\Lambda_{\alpha})}\int_{\Lambda_{\alpha}^{-1}%
\Lambda_{\alpha}}\left|  \gamma_{h}\right|  dh
\end{align*}
by Lemma 2.8.$~\square$

\bigskip

As opposed to Theorems 2.7 and 2.7$^{\prime}$, the measure in this proposition
has to be left invariant, hence when it is applied in tandem with Theorem 2.7
or 2.7$^{\prime}$, the measure will have to be invariant. Clearly Proposition
2.9 would also work if $\Lambda_{\alpha}$ wasn't necessarily open, but we had
$\Lambda_{\alpha}^{-1}\Lambda_{\alpha}\subset S_{\alpha}$ with $S_{\alpha}$
measurable and $\lim_{\alpha}\int_{S_{\alpha}}\left|  \gamma_{h}\right|  dh=0$.

\section{Weak mixing}

In this section we define weak mixing, and study some of its characterizations
using simple tools like density limits. This sets the stage for our study of
weak mixing of all orders in the next section. The discussion here is in a
fairly abstract setting, which for the most part does not require the net
$\left\{  \Lambda_{\alpha}\right\}  $ to be space-filling. As we will see, the
net is only required to be space-filling in order for the definition of weak
mixing to be independent of the net being used, and in the next section in the
final step of the proof of weak mixing to all orders, where the van der Corput
lemma is used.

\bigskip

\noindent\textbf{Definition 3.1. Dynamical system, measure preserving
dynamical system.} Let $(X,\Sigma,\nu)$ be a probability space. Let $K$ be any
semigroup. For each $g\in K$ let $T_{g}:X\rightarrow X$ be such that
$T_{g}\circ T_{h}=T_{gh}$ for all $g,h\in K$. Denote $g\mapsto T_{g}$ by $T$.
If $T_{g}^{-1}(\Sigma)\subset\Sigma$ for all $g\in K$, then $(X,\Sigma
,\nu,T,K)$ is called a \textit{dynamical system} (\textit{over }$K$; at times
it will be convenient to explicitly state the semigroup). If, additionally,
$\nu(T_{g}^{-1}(A))=\nu(A)$ for all $A\in\Sigma$ and $g\in K$, then
$(X,\Sigma,\nu,T,K)$ is called a \textit{measure preserving }dynamical system.

\bigskip

For a group (respectively semigroup) $G$, let $\mathfrak{M}_{G}$ denote the
set of all group (respectively semigroup) homomorphisms $G\rightarrow G$.

\bigskip

\noindent\textbf{Definition 3.2. Weak mixing and ergodicity.} Let $K$ be a
semigroup with a $\sigma$-algebra and a measure $\mu$. Let $\{\Lambda_{\alpha
}\}$ be a net of measurable subsets of $K$, such that $\mu(\Lambda_{\alpha
})>0$ for $\alpha$ large enough, and with $\mu(\Lambda_{\alpha})<\infty$ for
every $\alpha$. Let $M\subset\mathfrak{M}_{K}$. \noindent Assume that
$(X,\Sigma,\nu,T,K)$ is a dynamical system and that $g\mapsto\nu(A_{0}\cap
T_{\varphi(g)}^{-1}(A_{1}))$ is measurable for all $A_{0},A_{1}\in\Sigma$ and
all $\varphi\in M$.

\begin{description}
\item[(i)] $(X,\Sigma,\nu,T,K)$ is said to be $M$-\textit{weakly mixing}
\textit{relative t}o $\{\Lambda_{\alpha}\}$, if
\[
\lim_{\alpha}\frac{1}{\mu(\Lambda_{\alpha})}\int_{\Lambda_{\alpha}}\left|
\nu(A_{0}\cap T_{\varphi(g)}^{-1}(A_{1}))-\nu(A_{0})\nu(A_{1})\right|  dg=0
\]
for all $A_{0},A_{1}\in\Sigma$, and for all $\varphi\in M$.

\item[(ii)] $(X,\Sigma,\nu,T,K)$ is said to be $M$-\textit{ergodic}
\textit{relative to} $\{\Lambda_{\alpha}\}$, if
\[
\lim_{\alpha}\frac{1}{\mu(\Lambda_{\alpha})}\int_{\Lambda_{\alpha}}\nu
(A_{0}\cap T_{\varphi(g)}^{-1}(A_{1}))dg=\nu(A_{0})\nu(A_{1})
\]
for all $A_{0},A_{1}\in\Sigma$, and for all $\varphi\in M$.
\end{description}

\bigskip

\noindent\textbf{Remarks on Definition 3.2.} In the case of $K=\mathbb{N}$,
$\Lambda_{n}=\left\{  1,...,n\right\}  $ and with $M=\left\{  id_{\mathbb{N}%
}\right\}  $, Definition 3.2(i) corresponds to the usual definition of weak
mixing for an action of the semigroup $\mathbb{N}$, as given in (1.1). Since
all homomorphisms of $\mathbb{N}$ are of the form $n\mapsto kn$ for some
$k\in\mathbb{N}$, one can then easily show that $\left\{  id_{\mathbb{N}%
}\right\}  $-weak mixing implies $\mathfrak{M}_{\mathbb{N}}$-weak mixing.

For general $K$ our definition of weak mixing is quite abstract. We don't
assume the dynamical system to be measure preserving, or the net $\left\{
\Lambda_{\alpha}\right\}  $ to be space-filling in $K$, simply because these
assumptions are unnecessary in many of the results that follow, though they
are required when proving $M$-weak mixing of all orders. In ``practical''
cases that one usually studies in ergodic theory, one would expect these
assumptions to hold, for example $\Lambda_{n}$ mentioned above is
space-filling in $\mathbb{N}$. Under these assumptions, we will see in
Corollary 3.10 that the definition of weak mixing is independent of the
space-filling net we use, i.e. if a measure preserving dynamical system is $M
$-weakly mixing relative to one space-filling net, then it is $M$-weakly
mixing relative to all space-filling nets in $K$. The proof of Corollary 3.10,
as well as the parts of Propositions 3.8 and 3.9 which are used in this proof,
are the only places in this paper where we will use ergodicity.

In general the assumption that a dynamical system is $M$-weakly mixing, is a
restriction on $M$, since for example one would not expect $g\mapsto\nu
(A_{0}\cap T_{\varphi(g)}^{-1}(A_{1}))$ to\ even be measurable for all
homomorphisms $\varphi:K\rightarrow K$.

As a last remark, note that if $K$ has an identity $e$, and the homomorphism
given by $\varphi_{0}(g)=e$ for all $g\in G$ was in $M$, then the system
wouldn't be $M$-weakly mixing, hence we wouldn't want $\varphi_{0}$ to be in
$M$. We mention this simply because $\varphi_{0}$ does appear in the theory to
follow, but not as an element of $M$.

\bigskip

We now turn to a few technical tools which we will need in Section 4.

\bigskip

\noindent\textbf{Definition 3.3. Density zero, density limit.} Let $(G,\mu)$
be a measure space (with $G$ not necessarily a group or semigroup) and
$\{\Lambda_{\alpha}\}$ a net of measurable subsets of $G$. Assume that
$\mu(\Lambda_{\alpha})>0$ for $\alpha$ large enough, and that $\mu
(\Lambda_{\alpha})<\infty$ for every $\alpha$.

\begin{description}
\item[(i)] A set $R\subset G$ is said to have \textit{density zero relative
to} $\{\Lambda_{\alpha}\}$, and we write $D_{\{\Lambda_{\alpha}\}}(R)=0$ if
and only if there exists a measurable set $S\subset G$, with $R\subset S$ such
that
\[
\lim_{\alpha}\frac{\mu(\Lambda_{\alpha}\cap S)}{\mu(\Lambda_{\alpha})}=0\text{
\ \ .}%
\]

\item[(ii)] We say that $f:G\rightarrow L$, with $L$ a real or complex normed
space, has \textit{density limit} $a\in L$ relative to $\{\Lambda_{\alpha}\}$,
if and only if for each $\varepsilon>0$, $D_{\{\Lambda_{\alpha}\}}%
(S_{\varepsilon})=0$, where
\[
S_{\varepsilon}:=\{h\in G:\Vert f(h)-a\Vert\geq\varepsilon\}\text{ \ \ ,}%
\]
and we write it as
\[
D_{\{\Lambda_{\alpha}\}}\text{-}\lim f=D_{\{\Lambda_{\alpha}\}}\text{-}%
\lim_{h}f(h)=a\text{ \ \ .}%
\]
\end{description}

Note that if $R$ and $S$ have density zero relative to $\{\Lambda_{\alpha}\}$
and $V\subset S$, then $R\cap S$, $R\cup S$ and $V$ also have density zero
relative to $\{\Lambda_{\alpha}\}$.

\bigskip

\noindent\textbf{Proposition 3.4.}\emph{\ Let $f,g:G\rightarrow L$ with
$\left(  G,\mu\right)  $ and }$L$ \emph{as in Definition 3.3, and assume
that}
\[
D_{\{\Lambda_{\alpha}\}}\text{-}\lim f=a~~\mbox{ and }  ~~D_{\{\Lambda
_{\alpha}\}}\text{-}\lim g=b\text{ \ \ .}%
\]
\emph{Then}
\[
D_{\{\Lambda_{\alpha}\}}\text{-}\lim(f+g)=a+b
\]
\emph{and}
\[
D_{\{\Lambda_{\alpha}\}}\text{-}\lim(\beta f)=\beta a
\]
\emph{for any $\beta\in\mathbb{C}$. Furthermore, if $f,g$ are real-valued
functions and $f(h)\leq g(h)$ for all $h\in G$, then $a\leq b$.}

\bigskip

\noindent\textbf{Proof.} For each $\varepsilon>0$, let
\[
R_{\varepsilon}:=\{h\in G:\Vert f(h)-a\Vert\geq\varepsilon\}~~~\mbox{
and }
~~~S_{\varepsilon}:=\{h\in G:\Vert g(h)-b\Vert\geq\varepsilon\}\text{ \ \ .}%
\]
By definition, $R_{\varepsilon}$ and $S_{\varepsilon}$ have density zero
relative to $\{\Lambda_{\alpha}\}$. Let
\[
V_{\varepsilon}:=\{h\in G:\Vert(f+g)(h)-(a+b)\Vert\geq\varepsilon\}
\]
and
\[
V_{\varepsilon}^{\prime}:=\{h\in G:\Vert f(h)-a\Vert+\Vert g(h)-b\Vert
\geq\varepsilon\}\text{ \ \ .}%
\]
Since $\Vert(f+g)(h)-(a+b)\Vert\leq\Vert f(h)-a\Vert+\Vert g(h)-b\Vert$, it is
clear that $V_{\varepsilon}\subset V_{\varepsilon}^{\prime}$. Also, clearly
$V_{\varepsilon}^{\prime}\subset R_{\frac{\varepsilon}{2}}\cup
S_{\frac{\varepsilon}{2}}$. But $R_{\frac{\varepsilon}{2}}\cup
S_{\frac{\varepsilon}{2}}$ has density zero relative to $\{\Lambda_{\alpha}%
\}$, and hence the same holds for $V_{\varepsilon}^{\prime}$ and then
$V_{\varepsilon} $. Hence
\[
D_{\{\Lambda_{\alpha}\}}\text{-}\lim(f+g)=a+b\text{ \ \ .}%
\]

Letting $W_{\varepsilon}:=\{h\in G:\Vert(\beta f)(h)-\beta a\Vert
\geq\varepsilon\}$, it is easily seen that $W_{\varepsilon}$ has density zero
relative to $\{\Lambda_{\alpha}\}$, hence
\[
D_{\{\Lambda_{\alpha}\}}\text{-}\lim(\beta f)=\beta a\text{ \ \ .}%
\]
Finally, suppose that $f,g$ are real-valued functions, i.e. $L=\mathbb{R}$,
and $f(h)\leq g(h)$ for all $h\in G$. From the previous two results in this
proposition, we have that
\[
D_{\{\Lambda_{\alpha}\}}\text{-}\lim(g-f)=b-a\text{ \ \ .}%
\]
Hence for any $\varepsilon>0$, the set
\[
W_{\varepsilon}^{\prime}:=\{h\in G:|(g-f)(h)-(b-a)|\geq\varepsilon\}
\]
has density zero relative to $\{\Lambda_{\alpha}\}$. Suppose now that
$b-a=:\rho<0$. Since $(g-f)(h)\geq0$ for all $h\in G$, we must have that the
set $W_{\left|  \rho\right|  /2}^{\prime}$ consists of all of $G$. Hence
\[
\frac{\mu(\Lambda_{\alpha}\cap W_{|\rho|/2}^{\prime})}{\mu(\Lambda_{\alpha}%
)}=\frac{\mu(\Lambda_{\alpha})}{\mu(\Lambda_{\alpha})}=1\text{ \ \ ,}%
\]
contradicting the stated fact that $W_{|\rho|/2}^{\prime}$ has density zero
relative to $\{\Lambda_{\alpha}\}$. Therefore $b-a\geq0$. $\square$

\bigskip

We now give a Koopman-von Neumann type lemma:

\bigskip

\noindent\textbf{Lemma 3.5.}\emph{\ Let $(G,\mu)$ be a measure space, and let
$\{\Lambda_{\alpha}\}$ be a net of measurable subsets of $G$. Assume that
$\mu(\Lambda_{\alpha})>0$ for $\alpha$ large enough, and that $\mu
(\Lambda_{\alpha})<\infty$ for every $\alpha$. Let $f:G\rightarrow
\lbrack0,\infty)$ be bounded and measurable. Then the following are equivalent:}

\begin{description}
\item[(1)] $D_{\{\Lambda_{\alpha}\}}$-$\lim f=0$

\item[(2)] $\displaystyle        \lim_{\alpha}\frac{1}{\mu(\Lambda_{\alpha}%
)}\int_{\Lambda_{\alpha}}f~d\mu=0$
\end{description}

\noindent\textbf{Proof.} For every $\varepsilon>0$, let $S_{\varepsilon
}:=\{h\in G:f(h)\geq\varepsilon\}$, which is a measurable set, since $f$ is measurable.

(1) $\Rightarrow$ (2): From (1) we have that each $S_{\varepsilon}$ has
density zero relative to $\{\Lambda_{\alpha}\}$. Given any $\varepsilon>0$ and
index $\alpha$, consider the term
\[
\frac{1}{\mu(\Lambda_{\alpha})}\int_{\Lambda_{\alpha}}fd\mu=\frac{1}%
{\mu(\Lambda_{\alpha})}\int_{\Lambda_{\alpha}\cap S_{\varepsilon}}%
f~d\mu+\frac{1}{\mu(\Lambda_{\alpha})}\int_{\Lambda_{\alpha}\cap
S_{\varepsilon}^{c}}fd\mu\text{ \ \ .}%
\]
Since $S_{\varepsilon}$ has density zero relative to $\{\Lambda_{\alpha}\}$
\[
0\leq\frac{1}{\mu(\Lambda_{\alpha})}\int_{\Lambda_{\alpha}\cap S_{\varepsilon
}}f~dh\mu\leq\frac{\mu\left(  \Lambda_{\alpha}\cap S_{\varepsilon}\right)
}{\mu(\Lambda_{\alpha})}\sup f(G)\rightarrow0
\]
in the $\alpha$ limit. Also,
\[
0\leq\frac{1}{\mu(\Lambda_{\alpha})}\int_{\Lambda_{\alpha}\cap S_{\varepsilon
}^{c}}f~d\mu\leq\frac{\mu\left(  \Lambda_{\alpha}\cap S_{\varepsilon}%
^{c}\right)  }{\mu(\Lambda_{\alpha})}\varepsilon\leq\varepsilon\newline
\]
hence
\[
\lim_{\alpha}\frac{1}{\mu(\Lambda_{\alpha})}\int_{\Lambda_{\alpha}}%
fd\mu=0\text{ \ \ .}%
\]

(2) $\Rightarrow$ (1): Clearly $\varepsilon\chi_{S_{\varepsilon}}\leq f$. Also
note that $D_{\left\{  \Lambda_{\alpha}\right\}  }\left(  S_{\varepsilon
}\right)  =0$, since $S_{\varepsilon}$ is measurable and
\[
\varepsilon\frac{\mu\left(  \Lambda_{\alpha}\cap S_{\varepsilon}\right)  }%
{\mu\left(  \Lambda_{\alpha}\right)  }\leq\frac{1}{\mu(\Lambda_{\alpha})}%
\int_{\Lambda_{\alpha}}fd\mu
\]
which tends to zero in the $\alpha$ limit.$~\square$

\bigskip

\noindent\textbf{Corollary 3.6.} \emph{Consider the situation in Lemma 3.5,
except that we use $f:G\rightarrow\mathbb{R}$, assumed to be bounded and
measurable. Then}
\[
\lim_{\alpha}\frac{1}{\mu(\Lambda_{\alpha})}\int_{\Lambda_{\alpha}}%
[f(h)]^{2}~dh=0
\]
\emph{if and only if}
\[
\lim_{\alpha}\frac{1}{\mu(\Lambda_{\alpha})}\int_{\Lambda_{\alpha}%
}|f(h)|~dh=0\text{ \ \ .}%
\]

\bigskip

\noindent\textbf{Proof.} Given any $\varepsilon>0$. Let
\[
S_{\varepsilon}:=\{h\in G:[f(h)]^{2}\geq\varepsilon^{2}\}=\{h\in
G:|f(h)|\geq\varepsilon\}.
\]
Suppose that $\lim_{\alpha}\frac{1}{\mu(\Lambda_{\alpha})}\int_{\Lambda
_{\alpha}}[f(h)]^{2}~dh=0$, i.e. $D_{\{\Lambda_{\alpha}\}}$-$\lim
_{h}[f(h)]^{2}=0$ by Lemma 3.5. \noindent By the definition of the density
limit we have $D_{\{\Lambda_{\alpha}\}}(S_{\varepsilon})=0$. Since
$\varepsilon>0$ is arbitrary, we conclude that $D_{\{\Lambda_{\alpha}\}}%
$-$\lim|f|=0$, and hence $\lim_{\alpha}\frac{1}{\mu(\Lambda_{\alpha})}%
\int_{\Lambda_{\alpha}}|f(h)|~dh=0$ by Lemma 3.5.

The converse follows similarly.$~\square$

\bigskip

As a result, the $\left|  \cdot\right|  $ in Definition 3.2(i) of weak mixing,
can be replaced by $\left[  \cdot\right]  ^{2}$.

\bigskip

\noindent\textbf{Lemma 3.7.}\emph{\ Consider the situation in Lemma 3.5,
except that we use $f:G\rightarrow\mathbb{C}$, assumed to be bounded and
measurable. Let $\beta\in\mathbb{C}$.}

\bigskip\noindent If
\[
\lim_{\alpha}\frac{1}{\mu(\Lambda_{\alpha})}\int_{\Lambda_{\alpha}%
}f(h)dh=\beta
\]
\emph{and}
\[
\lim_{\alpha}\frac{1}{\mu(\Lambda_{\alpha})}\int_{\Lambda_{\alpha}}%
[f(h)]^{2}dh=\beta^{2}\text{ \ \ ,}%
\]
\emph{then}
\[
\lim_{\alpha}\frac{1}{\mu(\Lambda_{\alpha})}\int_{\Lambda_{\alpha}}\left[
f(h)-\beta\right]  ^{2}dh=0\text{ \ \ .}%
\]

\noindent\textbf{Proof.} This follows immediately if we note that
\begin{align*}
&  \frac{1}{\mu(\Lambda_{\alpha})}\int_{\Lambda_{\alpha}}\left[
f(h)-\beta\right]  ^{2}dh\\
&  =\frac{1}{\mu(\Lambda_{\alpha})}\int_{\Lambda_{\alpha}}\left(
[f(h)]^{2}-2\beta f(h)+\beta^{2}\right)  dh\\
&  \rightarrow0
\end{align*}
in the $\alpha$ limit.$~\square$

\bigskip

Next we consider standard characterizations of weak mixing, that we will need.
The first proposition does not require the system to be measure preserving,
but the second does.

\bigskip

\noindent\textbf{Proposition 3.8.} \emph{Let $K$ be a semigroup with a
$\sigma$-algebra and a measure $\mu$, and let $\{\Lambda_{\alpha}\}$ be a net
of measurable subsets of $K$. Assume that $\mu(\Lambda_{\alpha})>0$ for
$\alpha$ large enough, and that $\mu(\Lambda_{\alpha})<\infty$ for every
$\alpha$. Let $M\subset\mathfrak{M}_{K}$. Let $(X,\Sigma,\nu,T,K)$ be a
dynamical system. Set $(T\times T)_{h}=T_{h}\times T_{h}$ for all }$h\in
K$\emph{, where $(T_{h}\times T_{h})(x_{1},x_{2})=\left(  T_{h}(x_{1}%
),T_{h}(x_{2})\right)  $ for all $(x_{1},x_{2})\in X\times X$. Consider the
following statements:}

\begin{description}
\item[(1)] $(X,\Sigma,\nu,T,K)$\emph{\ is $M$-weakly mixing relative to
}$\{\Lambda_{\alpha}\}$.

\item[(2)] \emph{$(X\times X,\Sigma\times\Sigma,\nu\times\nu,T\times T,K)$ is
$M$-weakly mixing relative to $\{\Lambda_{\alpha}\}$.}

\item[(3)] \emph{$(X\times X,\Sigma\times\Sigma,\nu\times\nu,T\times T,K)$ is
$M$-ergodic relative to $\{\Lambda_{\alpha}\}$.}

\item[(4)] $\displaystyle   D_{\{\Lambda_{\alpha}\}}$-$\lim_{h}\nu(A_{0}\cap
T_{\phi(h)}^{-1}(A_{1}))=\nu(A_{0})\nu(A_{1})$\emph{\ for all }$A_{0},A_{1}%
\in\Sigma$\emph{\ and for each }$\varphi\in M$.
\end{description}

\noindent\emph{Then (1) and (4) are equivalent. Also, (2) implies (3), which
in turn implies(1).}

\bigskip

\noindent\textbf{Proof.} (1) $\Leftrightarrow$ (4): Given any $\varphi\in M$,
let $f(h):=\left|  \nu(A_{0}\cap T_{\varphi(h)}^{-1}(A_{1}))-\nu(A_{0}%
)\nu(A_{1})\right|  $, and apply Lemma 3.5.

\bigskip

\noindent(2) $\Rightarrow$ (3): Follows immediately from Definition 3.2.

\bigskip

\noindent(3) $\Rightarrow$ (1): Let $A_{0},A_{1}\in\Sigma$ and $\varphi\in M$.
We have
\begin{align*}
&  \lim_{\alpha}\frac{1}{\mu(\Lambda_{\alpha})}\int_{\Lambda_{\alpha}}%
\nu(A_{0}\cap T_{\varphi(g)}^{-1}(A_{1}))dg\\
&  =\lim_{\alpha}\frac{1}{\mu(\Lambda_{\alpha})}\int_{\Lambda_{\alpha}}%
(\nu\times\nu)((A_{0}\times X)\cap(T_{\varphi(g)}\times T_{\varphi(g)}%
)^{-1}(A_{1}\times X))dg\\
&  =(\nu\times\nu)(A_{0}\times X)(\nu\times\nu)(A_{1}\times X)\\
&  =\nu(A_{0})\nu(A_{1})\text{ \ \ ,}%
\end{align*}

\noindent and also%
\begin{align*}
&  \lim_{\alpha}\frac{1}{\mu(\Lambda_{\alpha})}\int_{\Lambda_{\alpha}}%
\nu(A_{0}\cap T_{\varphi(g)}^{-1}(A_{1}))^{2}dg\\
&  =\lim_{\alpha}\frac{1}{\mu(\Lambda_{\alpha})}\int_{\Lambda_{\alpha}}%
(\nu\times\nu)((A_{0}\times A_{0})\cap(T_{\varphi(g)}\times T_{\varphi
(g)})^{-1}(A_{1}\times A_{1}))dg\\
&  =(\nu\times\nu)(A_{0}\times A_{0})(\nu\times\nu)(A_{1}\times A_{1})\\
&  =\nu(A_{0})^{2}\nu(A_{1})^{2}\text{ \ \ .}%
\end{align*}
\noindent Therefore by Lemma 3.7 we have that
\[
\lim_{\alpha}\frac{1}{\mu(\Lambda_{\alpha})}\int_{\Lambda_{\alpha}}(\nu
(A_{0}\cap T_{\varphi(g)}^{-1}(A_{1}))-\nu(A_{0})\nu(A_{1}))^{2}dg=0,
\]
\noindent and it follows from Corollary 3.6 that $(X,\Sigma,\nu,T,K)$ is $M
$-weakly mixing relative to $\{\Lambda_{\alpha}\}$.$~\square$

\bigskip

\noindent\textbf{Proposition 3.9.} \emph{\ Consider the situation in
Proposition 3.8, but also assume that the dynamical system $(X,\Sigma
,\nu,T,K)$ is measure preserving. Then the following are equivalent:}

\begin{description}
\item[(1)] \emph{$(X,\Sigma,\nu,T,K)$ is $M$-weakly mixing relative to
$\{\Lambda_{\alpha}\}$.}

\item[(2)] \emph{$(X\times X,\Sigma\times\Sigma,\nu\times\nu,T\times T,K)$ is
$M$-weakly mixing relative to $\{\Lambda_{\alpha}\}$.}

\item[(3)] \emph{$(X\times X,\Sigma\times\Sigma,\nu\times\nu,T\times T,K)$ is
$M$-ergodic relative to $\{\Lambda_{\alpha}\}$.}

\item[(4)] $\displaystyle  \lim_{\alpha}\frac{1}{\mu(\Lambda_{\alpha})}%
\int_{\Lambda_{\alpha}}\left|  \langle f_{1},f_{2}\circ T_{\varphi(h)}%
\rangle-\langle f_{1},1\rangle\langle1,f_{2}\rangle\right|  dh=0$\emph{\ and
}$h\mapsto\langle f_{1},f_{2}\circ T_{\varphi(h)}\rangle$\emph{\ is
measurable\ for all }$f_{1},f_{2}\in L^{2}(\nu)$\emph{\ and for each }%
$\varphi\in M$.
\end{description}

\noindent\textbf{Proof.} By Proposition 3.8, we already have (2) $\Rightarrow$
(3) $\Rightarrow$ (1). Now for the rest:

\bigskip

\noindent(1) $\Rightarrow$ (2): Given any $\varphi\in M$ and $A,B,C,D\in
\Sigma$, we have
\begin{align*}
&  |\nu\times\nu((A\times C)\cap(T\times T)_{\varphi(h)}^{-1}(B\times
D))-\nu\times\nu(A\times C)\nu\times\nu(B\times D)|\\
&  =|\nu(A\cap T_{\varphi(h)}^{-1}(B))\nu(C\cap T_{\varphi(h)}^{-1}%
(D))-\nu(A)\nu(B)\nu(C)\nu(D)|\\
&  \leq\nu(A\cap T_{\varphi(h)}^{-1}(B))|\nu(C\cap T_{\varphi(h)}^{-1}%
(D))-\nu(C)\nu(D)|\\
&  +\nu(C)\nu(D)|\nu(A\cap T_{\varphi(h)}^{-1}(B))-\nu(A)\nu(B)|~\\
&  \leq\nu(A)|\nu(C\cap T_{\varphi(h)}^{-1}(D))-\nu(C)\nu(D)|\\
&  +\nu(C)\nu(D)|\nu(A\cap T_{\varphi(h)}^{-1}(B))-\nu(A)\nu(B)|\text{ \ \ .}%
\end{align*}
Hence by 3.8(1 and 4) and Proposition 3.4,%
\[
D_{\left\{  \Lambda_{\alpha}\right\}  }\text{-}\lim_{h}|\nu\times\nu((A\times
C)\cap(T\times T)_{\varphi(h)}^{-1}(B\times D))-\nu\times\nu(A\times
C)\nu\times\nu(B\times D)|=0\text{ \ \ .}%
\]
So again by 3.8(1 and 4), and since the system is measure preserving and the
rectangles form a semi-algebra that generates $\Sigma\times\Sigma$, the proof
follows in a standard way (see e.g. [12]).

\bigskip

\noindent\noindent(1) $\Rightarrow$ (4): This is true if $f_{1},f_{2}$ are
characteristic functions of measurable sets and given any $\varphi\in M$. The
desired result is obtained by forming linear combinations and approximating in
a standard way (see e.g. [12]).

\bigskip

\noindent(4) $\Rightarrow$ (1): This follows by taking $f_{1}$ and $f_{2}$ to
be characteristic functions of measurable sets, given any $\varphi\in
M$.$~\square$

\bigskip

We can now show that the definition of $M$-weak mixing relative to a
space-filling net, is independent of the space-filling net being used:

\bigskip

\noindent\textbf{Corollary 3.10.} \emph{\ If a measure preserving dynamical
system $(X,\Sigma,\nu,T,K)$ is $M$-weakly mixing relative to some
space-filling net in }$K$\emph{, then it is $M$-weakly mixing relative to
every space-filling net in }$K$\emph{.}

\bigskip

\noindent\textbf{Proof.} In [3] it is shown that the ergodicity of a
measure-preserving dynamical system is independent of the space-filling net
being used, and the proof holds for $M$-ergodicity, as we defined it here, as
well. Hence $M$-weak mixing is also independent of the space-filling net, by
the equivalence in Proposition 3.9(1 and 3).$~\square$

\section{Weak mixing of all orders}

In this section we show that weak mixing implies weak mixing of all orders.
Our approach is strongly influenced by that of [8] for the case of the group
$\mathbb{Z}$. The proof is by induction, two steps of which are given by the following:

\bigskip

\noindent\textbf{Proposition 4.1.}\emph{\ Let $K$ be a semigroup with a
$\sigma$-algebra and measure $\mu$, and assume that $K$ has an identity
element $e$. Let $\{\Lambda_{\alpha}\}$ be a net of measurable subsets of $K$
such that $\mu(\Lambda_{\alpha})>0$ for $\alpha$ large enough, and with
$\mu(\Lambda_{\alpha})<\infty$ for every $\alpha$. Let $M\subset
\mathfrak{M}_{K}$. Now we will use the following notation: $(X,\Sigma
,\nu,T,K)$ will denote any measure preserving dynamical system, but with }%
$K$\emph{\ fixed, and }$id$\emph{\ will denote the identity mapping
}$X\rightarrow X$\emph{. Let $\omega(f):=\int_{X}fd\nu$ for all $f\in
L^{\infty}(\nu)$, and let }$\left(  \omega\otimes\omega\right)  (f):=\int
_{X\times X}fd(\nu\times\nu)$\emph{\ for all $f\in L^{\infty}(\nu\times\nu)$.
Given $k\in\mathbb{N}$, let $\varphi_{1},...\varphi_{k}$ denote elements of
$M$, and let $f_{0},...,f_{k} $ denote real-valued elements of $L^{\infty}%
(\nu)$. Let $\varphi_{0}(h)=e$ for all $h\in K$.}

\emph{Consider the following statements (where the existence of the integrals
contained in each statement form part of that statement):}

\bigskip

\emph{1[}$\emph{k}$\emph{]: The integral $\int_{\Lambda_{\alpha}}\omega\left(
\prod_{j=0}^{k}f_{j}\circ T_{\varphi_{j}(g)}\right)  dg$ exists for all
$\alpha\geq\alpha_{0}$ for some $\alpha_{0}$, and}%
\[
\lim_{\alpha}\frac{1}{\mu(\Lambda_{\alpha})}\int_{\Lambda_{\alpha}}\left(
\omega\left(  \prod_{j=0}^{k}f_{j}\circ T_{\varphi_{j}(g)}\right)
-\prod_{j=0}^{k}\omega(f_{j})\right)  ^{2}dg=0\text{ \ \ .}%
\]

\emph{2[}$\emph{k}$\emph{]: }$\displaystyle   \lim_{\alpha}\frac{1}%
{\mu(\Lambda_{\alpha})}\int_{\Lambda_{\alpha}}\omega\left(  \prod_{j=0}%
^{k}f_{j}\circ T_{\varphi_{j}(g)}\right)  dg=\prod_{j=0}^{k}\omega(f_{j})$ \ \ .

\bigskip

\emph{3[}$\emph{k}$\emph{]: For $\kappa:=\prod_{j=1}^{k}\omega(f_{j})$, we have}%

\[
\lim_{\alpha}\left\|  \frac{1}{\mu(\Lambda_{\alpha})}\int_{\Lambda_{\alpha}%
}\prod_{j=1}^{k}f_{j}\circ T_{\varphi_{j}(g)}dg-\kappa\right\|  _{L^{2}(\nu
)}=0\text{ \ \ .}%
\]

\bigskip

\noindent\emph{Then}

\begin{description}
\item[(i)] \emph{1[}$\emph{k}$\emph{] implies 2[}$\emph{k}$\emph{].}

\item[(ii)] \emph{If 3[}$\emph{k}$\emph{] holds for all measure preserving
dynamical systems over }$K$\emph{\ with }$T_{e}=id$\emph{\ which are
$M$-weakly mixing relative to the given net $\left\{  \Lambda_{\alpha
}\right\}  $, and all $f_{1},...,f_{k}$ and all $\varphi_{1},...,\varphi_{k}$
with $\varphi_{j}\neq\varphi_{l}$ when $j\neq l$ for $j,l\in\left\{
1,...,k\right\}  $, then 1[}$\emph{k}$\emph{] also holds for all measure
preserving dynamical systems over }$K$\emph{\ with $T_{e}=id$ which are
$M$-weakly mixing relative to $\left\{  \Lambda_{\alpha}\right\}  $ and all
$f_{0},...,f_{k}$ and all $\varphi_{1},...,\varphi_{k}$ with $\varphi_{j}%
\neq\varphi_{l}$ when $j\neq l$ for $j,l\in\left\{  1,...,k\right\}  $.}
\end{description}

\noindent\textbf{Proof.} \noindent(i) Use Corollary 3.6.

(ii) The strong convergence in 3[$k$] implies weak convergence, i.e.
\begin{align*}
\lim_{\alpha}\left\langle \frac{1}{\mu(\Lambda_{\alpha})}\int_{\Lambda
_{\alpha}}\prod_{j=1}^{k}f_{j}\circ T_{\varphi_{j}(g)}~dg,~f_{0}\right\rangle
&  =\left\langle \kappa\cdot1,~f_{0}\right\rangle \\
&  =\omega\left(  \kappa f_{0}\right) \\
&  =\prod_{j=0}^{k}\omega(f_{j}).
\end{align*}
Furthermore, by the definition of the integral, and from the assumption that
$T_{\varphi_{0}(h)}=T_{e}=id$, we have that
\begin{align*}
&  \lim_{\alpha}\left\langle \frac{1}{\mu(\Lambda_{\alpha})}\int
_{\Lambda_{\alpha}}\prod_{j=1}^{k}f_{j}\circ T_{\varphi_{j}(g)}~dg,~f_{0}%
\right\rangle \\
&  =\lim_{\alpha}\frac{1}{\mu(\Lambda_{\alpha})}\left\langle \int
_{\Lambda_{\alpha}}\prod_{j=1}^{k}f_{j}\circ T_{\varphi_{j}(g)}~dg,~f_{0}%
\right\rangle \\
&  =\lim_{\alpha}\frac{1}{\mu(\Lambda_{\alpha})}\int_{\Lambda_{\alpha}%
}\left\langle \prod_{j=1}^{k}f_{j}\circ T_{\varphi_{j}(g)},~f_{0}\right\rangle
dg\\
&  =\lim_{\alpha}\frac{1}{\mu(\Lambda_{\alpha})}\int_{\Lambda_{\alpha}}%
\omega\left(  \prod_{j=0}^{k}f_{j}\circ T_{\varphi_{j}(g)}\right)  ~dg\text{
\ \ ,}%
\end{align*}
hence
\begin{equation}
\lim_{\alpha}\frac{1}{\mu(\Lambda_{\alpha})}\int_{\Lambda_{\alpha}}%
\omega\left(  \prod_{j=0}^{k}f_{j}\circ T_{\varphi_{j}(g)}\right)
~dg=\prod_{j=0}^{k}\omega(f_{j}) \tag{4.1.1}%
\end{equation}
and in particular the integral on the left exists for all $\alpha\geq
\alpha_{0}$ for some $\alpha_{0}$. Since by Proposition 3.9(1 and 2) the
product system $(X\times X,\Sigma\times\Sigma,\nu\times\nu,T\times T,K)$ is an
$M$-weak mixing dynamical system relative to $\{\Lambda_{\alpha}\}$, and the
product system is measure preserving with $(T\times T)_{e}=id\times id$, we
can apply (4.1.1) to the product system to obtain
\begin{align*}
&  \lim_{\alpha}\frac{1}{\mu(\Lambda_{\alpha})}\int_{\Lambda_{\alpha}}%
(\omega\otimes\omega)\left(  \prod_{j=0}^{k}\left(  f_{j}\otimes f_{j}\right)
\circ(T\times T)_{\varphi_{j}(g)}\right)  dg\\
&  =\prod_{j=0}^{k}(\omega\otimes\omega)(f_{j}\otimes f_{j}).
\end{align*}
where for every $f_{1},f_{2}\in L^{\infty}(\nu)$ we define $f_{1}\otimes
f_{2}:X\times X\rightarrow\mathbb{R}$ by $(f_{1}\otimes f_{2})(x_{1}%
,x_{2}):=f_{1}(x_{1})f_{2}(x_{2})$ for all $(x_{1},x_{2})\in X\times X$. By
Fubini's theorem, namely $(\omega\otimes\omega)\left(  f_{1}\otimes
f_{2}\right)  =\omega\left(  f_{1}\right)  \omega\left(  f_{2}\right)  $ for
all $f_{1},f_{2}\in L^{\infty}(\nu)$, we have
\[
\lim_{\alpha}\frac{1}{\mu(\Lambda_{\alpha})}\int_{\Lambda_{\alpha}}%
\omega\left(  \prod_{j=0}^{k}f_{j}\circ T_{\varphi_{j}(g)}\right)
^{2}dg=\prod_{j=0}^{k}\omega(f_{j})^{2}\text{ \ \ ,}%
\]
proving 1[$k$] by Lemma 3.7.~$\square$

\bigskip

Note that the only property of weak mixing which is used in Proposition 4.1's
proof, is that if a dynamical system is $M$-weakly mixing and measure
preserving, then so is its product with itself. This is the only reason that
the systems in Proposition 4.1 are required to be measure preserving,
otherwise Proposition 3.9(1 and 2) would not apply. Proposition 4.1 would
still hold if we considered dynamical systems with some abstract property,
call it $E$, instead of ``$M$-weak mixing and measure preserving'', as long as
the product of an $E$ dynamical system with itself is again an $E$ dynamical
system. In particular, even though Proposition 4.1 is expressed in terms of
functions instead of sets, we did not need the characterization of $M$-weak
mixing in terms of functions, given by Proposition 3.9(4).

In order to complete the induction argument, we need 1[$1$], and that if
2[$k-1$] holds for all measure preserving dynamical systems over $K$ with
$T_{e}=id$ which are $M$-weakly mixing relative to $\left\{  \Lambda_{\alpha
}\right\}  $, then the same is true for 3[$k$]. The latter requires some more
work, and we will need to specialize the $M$ that we will allow. Firstly note
that for an abelian group $G$ and any homomorphisms $\varphi_{1}$ and
$\varphi_{2}$ of $G$, the function $\varphi^{\prime}:G\rightarrow G$ defined
by
\begin{equation}
\varphi^{\prime}(g):=\varphi_{2}(g)\varphi_{1}(g)^{-1} \tag{4.1}%
\end{equation}
is also a homomorphism of $G$. Even though from now on we will use only
abelian groups, we will continue to use multiplicative notation, as in (4.1).

\bigskip

\noindent\textbf{Definition 4.2.} Let $G$ be an abelian group and let
$M\subset\mathfrak{M}_{G}$. We call $M$ \textit{translational} if for all
$\varphi_{1},\varphi_{2}\in M$ with $\varphi_{1}\neq\varphi_{2}$, the
homomorphism $\varphi^{\prime}$ defined by (4.1) is also in $M$.

\bigskip

\noindent\textbf{Proposition 4.3.} \emph{\ Let $G$ be an abelian group with a
$\sigma$-algebra and measure $\mu$, and let $M\subset\mathfrak{M}_{G}$ be
translational. Let $(X,\Sigma,\nu,T,G)$ be a measure preserving dynamical
system. Let $\omega(f):=\int_{X}fd\nu$ for all $f\in L^{\infty}(\nu)$. Let
$\left\{  \Lambda_{\alpha}\right\}  $ be a net of measurable subsets of $K$
with $\mu(\Lambda_{\alpha})<\infty$ for all $\alpha$, and $\mu(\Lambda_{\beta
})>0$ for $\beta$ large enough. Set $\varphi_{0}(g)=e$ for all $g\in G$.
Assume that for some $k\in\mathbb{N}$}
\[
\lim_{\alpha}\frac{1}{\mu(\Lambda_{\alpha})}\int_{\Lambda_{\alpha}}%
\omega\left(  \prod_{j=0}^{k-1}f_{j}\circ T_{\varphi_{j}(g)}\right)
dg=\prod_{j=0}^{k-1}\omega(f_{j})
\]
\emph{for all real-valued $f_{0},...,f_{k-1}\in L^{\infty}(\nu)$ and all
$\varphi_{1},...,\varphi_{k-1}\in M$ with $\varphi_{j}\neq\varphi_{l}$ when
$j\neq l$ for $j,l\in\left\{  1,...,k-1\right\}  $, and in particular the
existence of the integral over }$\Lambda_{\alpha}$\emph{\ and the limit is
assumed. Now set}
\[
u_{h}:=\prod_{j=1}^{k}f_{j}\circ T_{\varphi_{j}(h)}-\kappa
\]
\emph{for all $h\in K$, where $\kappa:=\prod_{j=1}^{k}\omega(f_{j})$, for a
given set of real-valued $f_{j}\in L^{\infty}(\nu)$ and $\varphi_{j}\in M$
with $\varphi_{j}\neq\varphi_{l}$ when $j\neq l$ for $j,l\in\left\{
1,...,k\right\}  $. Then}
\[
\gamma_{h}:=\lim_{\alpha}\frac{1}{\mu(\Lambda_{\alpha})}\int_{\Lambda_{\alpha
}}\left\langle u_{g},u_{gh}\right\rangle dg
\]
\emph{exists (where $\left\langle \cdot,\cdot\right\rangle $ is taken in
$L^{2}(\nu)$; $L^{\infty}(\nu)\subset L^{2}(\nu)$ since $\nu(X)<\infty$),
and}
\[
\gamma_{h}=\prod_{j=1}^{k}\omega\left(  f_{j}\left(  f_{j}\circ T_{\varphi
_{j}(h)}\right)  \right)  -\kappa^{2}%
\]
\emph{for all $h\in K$.}

\bigskip

\noindent\textbf{Proof.} We have%
\begin{align*}
&  \left\langle u_{g},u_{gh}\right\rangle \\
&  =\int_{X}\left(  \prod_{j=1}^{k}f_{j}\circ T_{\varphi_{j}(g)}%
-\kappa\right)  \left(  \prod_{l=1}^{k}f_{l}\circ T_{\varphi_{l}(gh)}%
-\kappa\right)  d\nu\\
&  =\int_{X}\left[  \prod_{j=1}^{k}\left(  f_{j}\circ T_{\varphi_{j}%
(g)}\right)  \left(  f_{j}\circ T_{\varphi_{j}(gh)}\right)  -\kappa\prod
_{j=1}^{k}f_{j}\circ T_{\varphi_{j}(g)}-\kappa\prod_{j=1}^{k}f_{j}\circ
T_{\varphi_{j}(gh)}\right]  d\nu+\kappa^{2}%
\end{align*}
and note that all three these last integrals exist, since $f_{j}\circ
T_{\varphi_{j}(g)}$ and products of such functions are in $L^{\infty}%
(\nu)\subset L^{1}(\nu)$. We now consider the three integrals in turn:

\bigskip

(a) Since $G$ is abelian and $T$ is measure preserving,
\begin{align*}
&  \int_{X}\left[  \prod_{j=1}^{k}\left(  f_{j}\circ T_{\varphi_{j}%
(g)}\right)  \left(  f_{j}\circ T_{\varphi_{j}(gh)}\right)  \right]  d\nu\\
&  =\int_{X}\left[  \prod_{j=1}^{k}\left(  f_{j}\circ T_{\varphi_{j}%
(g)}\right)  \left(  f_{j}\circ T_{\varphi_{j}(h)}\circ T_{\varphi_{j}%
(g)}\right)  \right]  d\nu\\
&  =\int_{X}\left\{  \prod_{j=1}^{k}\left[  f_{j}\left(  f_{j}\circ
T_{\varphi_{j}(h)}\right)  \right]  \circ T_{\varphi_{j}(g)\varphi_{1}%
(g)^{-1}}\right\}  \circ T_{\varphi_{1}(g)}d\nu\\
&  =\int_{X}\left\{  \prod_{j=1}^{k}\left[  f_{j}\left(  f_{j}\circ
T_{\varphi_{j}(h)}\right)  \right]  \circ T_{\varphi_{j}(g)\varphi_{1}%
(g)^{-1}}\right\}  d\left(  \nu\circ T_{\varphi_{1}(g)}^{-1}\right) \\
&  =\int_{X}\left\{  \prod_{j=1}^{k}\left[  f_{j}\left(  f_{j}\circ
T_{\varphi_{j}(h)}\right)  \right]  \circ T_{\varphi_{j}(g)\varphi_{1}%
(g)^{-1}}\right\}  d\nu\\
&  =\omega\left(  \prod_{j=0}^{k-1}\left[  f_{j+1}\left(  f_{j+1}\circ
T_{\varphi_{j+1}(h)}\right)  \right]  \circ T_{\varphi_{j}^{\prime}%
(g)}\right)
\end{align*}
where $\varphi_{j}^{\prime}(g):=\varphi_{j+1}(g)\varphi_{1}(g)^{-1}$ for all
$g\in G$ and $j=0,...,k-1$, so $\varphi_{j}^{\prime}\in M$ for $j=1,...,k-1$
since $M$ is translational, $\varphi_{j}^{\prime}\neq\varphi_{l}^{\prime}$
when $j\neq l$ for $j,l\in\left\{  1,...,k-1\right\}  $, and $\varphi
_{0}^{\prime}(g)=e$ for all $g\in G$. Hence%
\[
\lim_{\alpha}\frac{1}{\mu(\Lambda_{\alpha})}\int_{\Lambda_{\alpha}}\int
_{X}\left[  \prod_{j=1}^{k}\left(  f_{j}\circ T_{\varphi_{j}(g)}\right)
\left(  f_{j}\circ T_{\varphi_{j}(gh)}\right)  \right]  d\nu dg=\prod
_{j=0}^{k-1}\omega\left(  f_{j+1}\left(  f_{j+1}\circ T_{\varphi_{j+1}%
(h)}\right)  \right)
\]
by assumption.

\bigskip

(b) For the second integral, again using the fact that $T$ is measure
preserving, it follows as in (a) that%
\begin{align*}
\lim_{\alpha}\frac{1}{\mu(\Lambda_{\alpha})}\int_{\Lambda_{\alpha}}\int
_{X}\left[  \prod_{j=1}^{k}f_{j}\circ T_{\varphi_{j}(g)}\right]  d\nu dg  &
=\lim_{\alpha}\frac{1}{\mu(\Lambda_{\alpha})}\int_{\Lambda_{\alpha}}%
\omega\left(  \left[  \prod_{j=0}^{k-1}f_{j+1}\circ T_{\varphi_{j}^{\prime
}(g)}\right]  \right)  dg\\
&  =\prod_{j=0}^{k-1}\omega(f_{j+1})\\
&  =\kappa
\end{align*}
by assumption.

\bigskip

(c) Lastly, again since $G$ is abelian and $T$ is measure preserving,
\begin{align*}
&  \lim_{\alpha}\frac{1}{\mu(\Lambda_{\alpha})}\int_{\Lambda_{\alpha}}\int
_{X}\left[  \prod_{j=1}^{k}f_{j}\circ T_{\varphi_{j}(gh)}\right]  d\nu dg\\
&  =\lim_{\alpha}\frac{1}{\mu(\Lambda_{\alpha})}\int_{\Lambda_{\alpha}}%
\omega\left(  \prod_{j=0}^{k-1}\left(  f_{j+1}\circ T_{\varphi_{j+1}%
(h)}\right)  \circ T_{\varphi_{j}^{\prime}(g)}\right)  dg\\
&  =\prod_{j=0}^{k-1}\omega(f_{j+1}\circ T_{\varphi_{j+1}(h)})\\
&  =\prod_{j=0}^{k-1}\omega(f_{j+1})\\
&  =\kappa
\end{align*}
by assumption.

\bigskip

(d) From (a)-(c)%
\[
\gamma_{h}=\prod_{j=1}^{k}\omega\left(  f_{j}\left(  f_{j}\circ T_{\varphi
_{j}(h)}\right)  \right)  -\kappa^{2}%
\]
and in particular $\gamma_{h}$ exists.$~\square$

\bigskip

We now state and prove our final result, namely that weak mixing implies weak
mixing of all orders. This is where our van der Corput lemma is finally
applied, along with Propositions 4.1 and 4.3, and the characterization of
$M$-weak mixing given by Proposition 3.9(1 and 4) which so far we have not used.

\bigskip

\noindent\textbf{Theorem 4.4.} \emph{\ Let $(X,\Sigma,\nu,T,G)$ be a measure
preserving dynamical system for an abelian second countable topological group
$G$ with invariant measure $\mu$, and with $T_{e}=id$. Let $M\subset
\mathfrak{M}_{G}$ be translational. Assume that $(X,\Sigma,\nu,T,G)$ is
$M$-weakly mixing relative to a uniformly space-filling sequence of open sets
$\{\Lambda_{n}\}$ in $G$, and that $(X,\Sigma,\nu,T,G)$ is $M$-weakly mixing
relative to the sequence $\left\{  \Lambda_{n}^{-1}\Lambda_{n}\right\}  $, so
in particular we require }$\mu(\Lambda_{n}^{-1}\Lambda_{n})>0$\emph{\ for }$n
$\emph{\ large enough and }$\mu(\Lambda_{n}^{-1}\Lambda_{n})<\infty
$\emph{\ for every }$n$\emph{, and where we assume that $\mu(\Lambda_{n}%
^{-1}\Lambda_{n})\leq c\mu(\Lambda_{n})$ for $n$ large enough and some
strictly positive real number $c$. Assume furthermore that $G\rightarrow
L^{\infty}(\nu):g\mapsto f\circ T_{\varphi(g)}$ is continuous in the
$L^{\infty}$-norm topology on $L^{\infty}(\nu)$ for all $\varphi\in M$. Then}
\[
\lim_{n\rightarrow\infty}\frac{1}{\mu(\Lambda_{n})}\int_{\Lambda_{n}}\left(
\omega\left(  \prod_{j=0}^{k}f_{j}\circ T_{\varphi_{j}(g)}\right)
-\prod_{j=0}^{k}\omega(f_{j})\right)  ^{2}dg=0
\]
\emph{for any real-valued $f_{j}\in L^{\infty}(\nu)$ and any $\varphi
_{1},...,\varphi_{k}\in M$ with $\varphi_{j}\neq\varphi_{l}$ when $j\neq l$
for $j,l\in\left\{  1,...,k\right\}  $, and with $\varphi_{0}(g)=e$ for all
$g\in G$, where $\omega(f):=\int_{X}fd\nu$.}

\bigskip

\noindent\textbf{Proof.} We need to complete the induction argument started in
Proposition 4.1, and we will continue using its notation, but with $K=G$.
Since $G\rightarrow L^{\infty}(\nu):g\mapsto f\circ T_{\varphi(g)}$ is
continuous, so is $F:G\rightarrow L^{\infty}(\nu):g\mapsto\prod_{j=0}^{k}%
f_{j}\circ T_{\varphi_{j}(g)}$ in the $L^{\infty}$-topology. Since $\nu(X)=1$,
we have $\left|  \left|  f\right|  \right|  _{2}\leq\left|  \left|  f\right|
\right|  _{\infty}$ for all $f\in L^{\infty}(\nu)$, so the $L^{\infty}%
$-topology is finer that the $L^{2}$-topology, hence $F$ is continuous in the
$L^{2}$-topology on $L^{\infty}(\nu)$ as well. It follows that
\[
G\times G\rightarrow\mathbb{R}:(g,h)\mapsto\left\langle \prod_{j=1}^{k}%
f_{j}\circ T_{\varphi_{j}(g)},\prod_{j=1}^{k}f_{j}\circ T_{\varphi_{j}%
(h)}\right\rangle
\]
is continuous. Keep in mind that $\omega\left(  \left(  \prod_{j=1}^{k}%
f_{j}\circ T_{\varphi_{j}(g)}\right)  \left(  \prod_{j=1}^{k}f_{j}\circ
T_{\varphi_{j}(h)}\right)  \right)  =\left\langle \prod_{j=1}^{k}f_{j}\circ
T_{\varphi_{j}(g)},\prod_{j=1}^{k}f_{j}\circ T_{\varphi_{j}(h)}\right\rangle
$. Now we write
\[
u_{g}:=\prod_{j=1}^{k}f_{j}\circ T_{\varphi_{j}(g)}-\kappa
\]
for all $g\in G$, where $\kappa:=\prod_{j=1}^{k}\omega(f_{j})$. It follows
that $G\times G\rightarrow\mathbb{C}:(g,h)\mapsto\left\langle u_{g}%
,u_{h}\right\rangle $ is continuous and therefore Borel measurable. Note that
$g\mapsto\left\langle u_{g},x\right\rangle $ is also Borel measurable for all
$x\in L^{2}(\nu)$. Furthermore, $G\rightarrow L^{2}(\nu):g\mapsto u_{g}$ is
bounded, since each $f_{j}$ is essentially bounded and $\nu(X)=1$. (We need
these properties, since we will be applying Theorem 2.7$^{\prime}$ to the
function $g\mapsto u_{g}$.) Since $\mu(\Lambda_{n}^{-1}\Lambda_{n})\leq
c\mu(\Lambda_{n})$, and we have $M$-weak mixing relative to $\left\{
\Lambda_{n}^{-1}\Lambda_{n}\right\}  $, it follows from Proposition 3.9(1 and
4) that
\begin{equation}
\lim_{n\rightarrow\infty}\frac{1}{\mu(\Lambda_{n})}\int_{\Lambda_{n}%
^{-1}\Lambda_{n}}\left|  \omega\left(  f_{0}\left(  f_{1}\circ T_{\varphi
(g)}\right)  \right)  -\omega(f_{0})\omega(f_{1})\right|  dg=0 \tag{4.4.1}%
\end{equation}
for all $\varphi\in M$. By Proposition 4.3, assuming 2[$k-1$] for all measure
preserving dynamical systems over $G$ with $T_{e}=id$, which are $M$-weakly
mixing relative to $\left\{  \Lambda_{\alpha}\right\}  $, and of course for
all $f_{0},...,f_{k-1}$ and all $\varphi_{1},...,\varphi_{k-1}\in M$ with
$\varphi_{j}\neq\varphi_{l}$ when $j\neq l$ for $j,l\in\left\{
1,...,k-1\right\}  $, we have%
\[
\gamma_{h}:=\lim_{n\rightarrow\infty}\frac{1}{\mu(\Lambda_{n})}\int
_{\Lambda_{n}}\left\langle u_{g},u_{gh}\right\rangle dg=\prod_{j=1}^{k}%
\omega\left(  f_{j}\left(  f_{j}\circ T_{\varphi_{j}(h)}\right)  \right)
-\prod_{j=1}^{k}\omega(f_{j})^{2}%
\]
for any $f_{1},...,f_{k}$ and all $\varphi_{1},...,\varphi_{k}\in M$ with
$\varphi_{j}\neq\varphi_{l}$ when $j\neq l$ for $j,l\in\left\{
1,...,k\right\}  $, for all $h\in G$. Using the identity $\prod_{j=1}^{k}%
a_{j}-\prod_{j=1}^{k}b_{j}=\sum_{j=1}^{k}\left(  \prod_{l=1}^{j-1}%
a_{l}\right)  \left(  a_{j}-b_{j}\right)  \left(  \prod_{l=j+1}^{k}%
b_{l}\right)  $ it follows that
\[
\int_{\Lambda_{m}^{-1}\Lambda_{m}}\left|  \gamma_{h}\right|  dh\leq\sum
_{j=1}^{k}A_{j}\left|  \prod_{l=j+1}^{k}\omega(f_{l})^{2}\right|
\int_{\Lambda_{m}^{-1}\Lambda_{m}}\left|  \omega\left(  f_{j}\left(
f_{j}\circ T_{\varphi_{j}(h)}\right)  \right)  -\omega(f_{j})^{2}\right|  dh
\]
where $A_{j}:=\sup_{h\in G}\left|  \prod_{l=1}^{j-1}\omega\left(  f_{l}\left(
f_{l}\circ T_{\varphi_{l}(h)}\right)  \right)  \right|  $ which exists in
$\mathbb{R}$, since the $f_{j}$ 's are essentially bounded. Note that
$\int_{\Lambda_{m}^{-1}\Lambda_{m}}\left|  \gamma_{h}\right|  dh$ exists,
since the integrand is continuous. Hence
\[
\lim_{m\rightarrow\infty}\frac{1}{\mu(\Lambda_{m})}\int_{\Lambda_{m}%
^{-1}\Lambda_{m}}\left|  \gamma_{h}\right|  dh=0
\]
by (4.4.1). From Proposition 2.9 and Theorem 2.7$^{\prime}$ we then have
\[
\lim_{n\rightarrow\infty}\frac{1}{\mu(\Lambda_{n)}}\int_{\Lambda_{n}}%
u_{g}dg=0
\]
where the limit is taken in the $L^{2}$-norm, i.e. 3[$k$] holds for all
measure preserving dynamical systems over $G$ with $T_{e}=id$, which are
$M$-weakly mixing relative to $\left\{  \Lambda_{\alpha}\right\}  $, and all
$f_{1},...,f_{k}$ and all $\varphi_{1},...,\varphi_{k}\in M$ with $\varphi
_{j}\neq\varphi_{l}$ when $j\neq l$ for $j,l\in\left\{  1,...,k\right\}  $.
But 1[$1$] holds for all $f_{0},f_{1}\in L^{\infty}(\nu)$ and all $\varphi\in
M$ for all measure preserving dynamical systems over $G$ with $T_{e}=id$,
which are $M$-weakly mixing relative to $\left\{  \Lambda_{\alpha}\right\}  $,
because of Proposition 3.9(1 and 4) and Corollary 3.6, completing the
induction argument started in Proposition 4.1, and proving 1[$k$] for all
$k\in\mathbb{N}$.$~\square$

\bigskip

By Corollary 3.6, the $\left[  \cdot\right]  ^{2}$ in the integrand in Theorem
4.4, can be replaced by $\left|  \cdot\right|  $, to have the same form as
Definition 3.2(i) of weak mixing.

Note that if $\left\{  \Lambda_{n}^{-1}\Lambda_{n}\right\}  $ is also
space-filling in $G$, then the assumption that the system be $M$-weakly mixing
relative to $\left\{  \Lambda_{n}^{-1}\Lambda_{n}\right\}  $ can be dropped
because of Corollary 3.10. If $\left\{  \Lambda_{n}^{-1}\Lambda_{n}\right\}  $
does not have the properties required in Theorem 4.4, for example if the
system is not $M$-weak mixing relative to $\left\{  \Lambda_{n}^{-1}%
\Lambda_{n}\right\}  $, but there is some other uniformly space-filling
sequence $\left\{  \Lambda_{n}^{\prime}\right\}  $ such that $\left\{
\Lambda_{n}^{\prime-1}\Lambda_{n}^{\prime}\right\}  $ does have the required
properties, then we can replace $\left\{  \Lambda_{n}\right\}  $ by $\left\{
\Lambda_{n}^{\prime}\right\}  $ because of Corollary 3.10, to get weak mixing
of all orders relative to $\left\{  \Lambda_{n}^{\prime}\right\}  $. We now
briefly consider examples of space-filling sequences with the required properties.

In the simple case where $G=\mathbb{Z}$ with the counting measure $\mu$, and
$\Lambda_{n}=\{-n,\ldots,n\}$ which is uniformly space-filling in $\mathbb{Z}
$, we have $\Lambda_{n}^{-1}\Lambda_{n}=\{-2n,\ldots,2n\}$, so $\mu
(\Lambda_{n})\leq\mu(\Lambda_{n}^{-1}\Lambda_{n})\leq2\mu(\Lambda_{n})$ for
$n\geq1$, and if the dynamical system is weak mixing relative to
$\{\Lambda_{n}\}$, then it is also weak mixing relative to $\Lambda_{n}%
^{-1}\Lambda_{n}=\Lambda_{2n}$. Hence the conditions of Theorem 4.4 are
satisfied. Furthermore, if the system is only weak mixing relative to
$\{0,\ldots,n\}$, so we are working over the semigroup $\mathbb{N}\cup\{0\}$,
and $T$ is injective, then it is easily seen that it is also weak mixing
relative to $\Lambda_{n}$. This implies the usual version of weak mixing of
all orders when working on the semigroup $\mathbb{N}\cup\{0\}$, for an
injective $T$.

As another example of a sequence with the properties in Theorem 4.4, let
$\Lambda_{m}$ be the open ball of radius $m$ in $\mathbb{R}^{q}$ for any
positive integer $q$. Note that $\left\{  \Lambda_{m}\right\}  $ is a
uniformly space-filling sequence in $\mathbb{R}^{q}$. Then $\Lambda_{m}%
^{-1}\Lambda_{m}=\Lambda_{2m}$, which means that $M$-weak mixing relative to
$\left\{  \Lambda_{m}\right\}  $, implies $M$-weak mixing relative to
$\left\{  \Lambda_{m}^{-1}\Lambda_{m}\right\}  $, while $\mu(\Lambda_{m}%
^{-1}\Lambda_{m})=2^{q}\mu(\Lambda_{m})$, as is required in Theorem 4.4.

Concerning the assumption that $M$ is translational, a simple example would be
of the following type: Use the group $G=\mathbb{R}^{q}$. Let $M$ be all
$q\times q$ non-zero diagonal real matrices acting as linear operators on
$\mathbb{R}^{q}$. (We exclude the zero matrix simply because this would make
$M$-weak mixing impossible.) Then $M$ is a translational set of homomorphisms
of $\mathbb{R}^{q}$. The same is true if we drop the condition that the
matrices be diagonal. Similarly if we work with $\mathbb{Z}^{q}$ instead of
$\mathbb{R}^{q}$ and use matrices over the integers.

\bigskip

\textbf{Acknowledgment.} We thank Richard de Beer, Willem Fouch\'{e}, Johan
Swart and Gusti van Zyl for useful discussions.

\end{document}